\documentclass[a4paper,12pt,reqno]{amsart}
\usepackage{amsmath}
\usepackage{amssymb}
\usepackage{amsthm}

\usepackage{amscd}
\usepackage{amsfonts}
\usepackage{setspace}
\usepackage{version}

\title[Bombieri--Vinogradov theorems for smooth numbers]{Bombieri--Vinogradov and Barban--Davenport--Halberstam type theorems for smooth numbers}
\author{Adam J Harper}
\address{Department of Pure Mathematics and Mathematical Statistics, Wilberforce Road, Cambridge CB\textup{3} \textup{0}WA, England}
\email{A.J.Harper@dpmms.cam.ac.uk}
\date{29th August 2012}
\thanks{The author is supported by a Doctoral Prize from the Engineering and Physical Sciences Research Council of the United Kingdom.}

\numberwithin{equation}{section}
\theoremstyle{plain}

\onehalfspace
\newcommand{\N}{\mathbb{N}}

\newtheorem{thm1}{Theorem}
\newtheorem{thm2}[thm1]{Theorem}
\newtheorem{thm3}[thm1]{Theorem}
\newtheorem{smooth1}{Smooth Numbers Result}
\newtheorem{smooth2}[smooth1]{Smooth Numbers Result}
\newtheorem{smooth3}[smooth1]{Smooth Numbers Result}
\newtheorem{zeros1}{Zeros Result}
\newtheorem{prop1}{Proposition}
\newtheorem{lem1}{Lemma}
\newtheorem{lem2}[lem1]{Lemma}
\newtheorem{multls}{Multiplicative Large Sieve}
\newtheorem{prop2}[prop1]{Proposition}

\begin{document}

\maketitle

\begin{abstract}
We prove Bombieri--Vinogradov and Barban--Davenport--Halberstam type theorems for the $y$-smooth numbers less than $x$, on the range $\log^{K}x \leq y \leq x$. This improves on the range $\exp\{\log^{2/3+\epsilon}x\} \leq y \leq x$ that was previously available. Our proofs combine zero-density methods with direct applications of the large sieve, which seems to be an essential feature and allows us to cope with the sparseness of the smooth numbers. We also obtain improved individual (i.e. not averaged) estimates for character sums over smooth numbers.
\end{abstract}

\section{Introduction}
For $y \geq 1$, let $\mathcal{S}(y)$ denote the set of $y$-smooth numbers: that is, the set of numbers all of whose prime factors are less than or equal to $y$. Smooth numbers are ubiquitous in analytic number theory, (see e.g. the survey paper~\cite{ht3} of Hildebrand and Tenenbaum), and it is natural to investigate many of the same questions for them as are studied for the prime numbers. For example, one might be interested in the distribution of smooth numbers among the integers less than $x$; in arithmetic progressions; in short intervals $[x,x+z]$; or in arithmetic progressions on average. In this paper we will prove some results concerning the latter problem.

In the case of primes, the most celebrated theorem concerning their average distribution in arithmetic progressions is undoubtedly the {\em Bombieri--Vinogradov theorem}: for any fixed $A > 0$, and any $\sqrt{x} \log^{-A}x \leq Q \leq \sqrt{x}$, we have
$$ \sum_{q \leq Q} \max_{y \leq x} \max_{(a,q)=1} \left| \Psi(y;q,a)-\frac{y}{\phi(q)} \right| \ll_{A} \sqrt{x}Q\log^{4}x, $$
where as usual we set $\Psi(y;q,a) := \sum_{n \leq y, n \equiv a (\textrm{mod } q)} \Lambda(n)$. In particular, if $Q = \sqrt{x} \log^{-A}x$ then the right hand side is $\ll_{A} x\log^{4-A}x$, which beats the trivial bound $\ll x\log Q$ (obtained since $\Psi(x;q,a) \ll x/\phi(q)$ for $q \leq \sqrt{x}$) provided $A > 3$.

If one is prepared to sum over residue classes $a$, rather than taking a maximum, then one can obtain an interesting bound with $Q$ much larger. Such a result may be called a {\em Barban--Davenport--Halberstam theorem}: for any fixed $A > 0$, and any $x \log^{-A}x \leq Q \leq x$, we have
$$ \sum_{q \leq Q} \sum_{(a,q)=1} \left| \Psi(x;q,a)-\frac{x}{\phi(q)} \right|^{2} \ll_{A} xQ \log x . $$

The theorems quoted above are essentially as stated in Davenport's book~\cite{davenport}, and are actually due to Vaughan~\cite{vaughan} and to Gallagher~\cite{gallagher}, respectively\footnote{In Davenport's book~\cite{davenport} the Bombieri--Vinogradov theorem is stated with a bound $\sqrt{x}Q\log^{5}x$, but the proof is by Vaughan's~\cite{vaughan} method which readily yields a bound $\sqrt{x}Q\log^{4}x$ if one is slightly more careful.}. The original results of Bombieri~\cite{bomb} and Vinogradov~\cite{aivino}, and of Barban (see $\S 3$ of \cite{barban}) and Davenport and Halberstam~\cite{davhal}, were slightly quantitatively weaker.

Let us point out two things about these results that will be significant later. Firstly, both bounds are ineffective (unless $A$ is very small), because the proofs rely on information about possible Siegel zeros of Dirichlet $L$-functions. In our current state of knowledge this ineffectiveness seems unavoidable, because if a Siegel zero did exist it would genuinely distort the distribution of primes in arithmetic progressions. Secondly, the known proofs actually give bounds for absolute values of character sums, namely
$$ \sum_{q \leq Q} \frac{1}{\phi(q)} \sum_{\substack{\chi \; (\textrm{mod } q), \\ \chi \neq \chi_{0}}} \left|\sum_{n \leq x} \Lambda(n) \chi(n) \right| \;\;\; \textrm{ and } \;\;\; \sum_{q \leq Q} \frac{1}{\phi(q)} \sum_{\substack{\chi \; (\textrm{mod } q), \\ \chi \neq \chi_{0}}} \left|\sum_{n \leq x} \Lambda(n) \chi(n) \right|^{2}, $$
respectively. One wouldn't hope for a better bound than $|\sum_{n \leq x} \Lambda(n) \chi(n)| \ll \sqrt{x}$, and the right hand sides in the theorems do indeed correspond, up to logarithmic factors, to such squareroot cancellation. Thus, although we would expect a non-trivial Bombieri--Vinogradov type theorem to hold with $Q$ much larger than $\sqrt{x}$, we wouldn't hope to prove such a result in this way (i.e. by bounding absolute values of character sums). In the case of the Barban--Davenport--Halberstam theorem the bound is known to be sharp: see e.g. the paper~\cite{mont2} of Montgomery. Indeed, in that case the original problem is equivalent to bounding the character sum that we wrote down.

\vspace{12pt}
Now we turn to the case of smooth numbers. For $x \geq 1$, and natural numbers $a,q$, we define
$$ \Psi_{q}(x,y) := \sum_{n \leq x, (n,q)=1} \textbf{1}_{\{n \in \mathcal{S}(y)\}} \;\;\; \textrm{ and } \;\;\; \Psi(x,y;q,a) := \sum_{n \leq x, n \equiv a (\textrm{mod } q)} \textbf{1}_{\{n \in \mathcal{S}(y)\}}, $$
where $\textbf{1}$ denotes the indicator function. We will also write $\Psi(x,y) := \sum_{n \leq x} \textbf{1}_{\{n \in \mathcal{S}(y)\}}$. Thus our task is to give bounds, on average, for $|\Psi(x,y;q,a) - \Psi_{q}(x,y)/\phi(q)|$. As in the case of prime numbers, one of our concerns will be to obtain good bounds with as large a range of summation over $q$ as possible. However, here we also need to be concerned about the range of $y$ that we can handle. When $y$ is small compared with $x$ the $y$-smooth numbers become a very sparse set, (for example if $y = \log^{K}x$ for some fixed $K \geq 1$ then $\Psi(x,\log^{K}x)/x \approx x^{-1/K}$, which is less by far than the density $\approx 1/\log x$ of the primes less than $x$), which creates new and interesting difficulties.

We shall prove the following results: in their statements $c > 0$ and $K > 0$ are certain fixed and effective constants, that should be thought of as quite small and quite large, respectively.
\begin{thm1}
Let $\log^{K}x \leq y \leq x$ be large and let $1 \leq Q \leq \sqrt{\Psi(x,y)}$. Then
$$ \sum_{q \leq Q} \max_{(a,q)=1} \left| \Psi(x,y;q,a)-\frac{\Psi_{q}(x,y)}{\phi(q)} \right| \ll \Psi(x,y)\left(e^{-\frac{cu}{\log^{2}(u+1)}} + y^{-c}\right) + \sqrt{\Psi(x,y)}Q\log^{7/2}x, $$
where $u:=(\log x)/\log y$, and the implicit constant is absolute {\em and effective}. In addition, for any $A > 0$ we have the bound
$$ \ll_{A} \Psi(x,y)\left(\frac{e^{-\frac{cu}{\log^{2}(u+1)}}}{\log^{A}x} + y^{-c}\right) + \sqrt{\Psi(x,y)}Q\log^{7/2}x, $$
but the implicit constant is now ineffective.
\end{thm1}

\begin{thm2}
Let $\log^{K}x \leq y \leq x$ be large and let $1 \leq Q \leq \Psi(x,y)$. Then
$$ \sum_{q \leq Q} \sum_{(a,q)=1} \left| \Psi(x,y;q,a)-\frac{\Psi_{q}(x,y)}{\phi(q)} \right|^{2} \ll \Psi(x,y)^{2} \left(e^{-\frac{cu}{\log^{2}(u+1)}} + y^{-c}\right) + \Psi(x,y)Q, $$
where the implicit constant is absolute {\em and effective}. In addition, for any $A > 0$ we have the bound
$$ \ll_{A} \Psi(x,y)^{2}\left(\frac{e^{-\frac{cu}{\log^{2}(u+1)}}}{\log^{A}x} + y^{-c}\right) + \Psi(x,y)Q , $$
but the implicit constant is now ineffective.
\end{thm2}

Theorems 1 and 2 supply non-trivial bounds when $Q=\sqrt{\Psi(x,y)}\log^{-A}x$ and $Q=\Psi(x,y)\log^{-A}x$, respectively. This mirrors the classical Bombieri--Vinogradov and Barban--Davenport--Halberstam theorems, and as there one couldn't hope to prove non-trivial bounds for larger $Q$ by a method based on absolute values of character sums\footnote{By ``non-trivial bounds'' we mean ``better than would follow if we just assumed that $|\Psi(x,y;q,a)-\Psi_{q}(x,y)/\phi(q)| \ll \Psi_{q}(x,y)/\phi(q)$ always''. However, for general $x,y,q$ we cannot prove bounds like this, (see Hildebrand and Tenenbaum's survey paper~\cite{ht3} for some discussion of what is known), so Theorems 1 and 2 might still be interesting for larger $Q$.}. Moreover, provided that $y$ isn't too big, or more specifically provided $u/(\log^{2}(u+1) \log\log x) \rightarrow \infty$, we save an arbitrary power of a logarithm {\em in an effective way}. We can prove strong effective results because a hypothetical Siegel zero does not distort the distribution of $y$-smooth numbers too much, when $y$ is small. One way to think about this is to note that if the $L$-series $L(s,\chi)$ has a Siegel zero, then the character $\chi$ must ``behave a lot like'' the M\"{o}bius function $\mu$. See e.g. Granville and Soundararajan's article~\cite{gransound} for some discussion of such issues. Whereas the sum of the M\"{o}bius function over primes is very large, the sum of the M\"{o}bius function over smooth numbers exhibits cancellation, as discussed in detail in Tenenbaum's paper~\cite{tenerdosall}.

Theorem 1 improves on Th\'{e}or\`{e}me 1 of Fouvry and Tenenbaum~\cite{fouvryten}, who proved such a result on the restricted range $\exp\{\log^{2/3+\epsilon}x\} \leq y \leq x$, for any fixed $\epsilon > 0$, with a bound roughly of the form
$$ \ll_{A,\epsilon} \frac{\Psi(x,y) e^{-c_{1}u/\log^{2}(u+1)}}{\log^{A}x} + \sqrt{x}Q u^{u(1+o_{A,\epsilon}(1))} \log^{A+5}x . $$
Here $c_{1}=c_{1}(\epsilon,A) > 0$ is a certain constant, and the $o(1)$ term tends to 0 (in a way that depends a little bit on $\epsilon$ and $A$) as $u \rightarrow \infty$. 
Note that the bound in Theorem 1 is always at least as good as this. Although Fouvry and Tenenbaum do not emphasise it, the implicit constant in their bound appears, like ours, to be effectively computable when $A$ is small.

See Fouvry and Tenenbaum's paper~\cite{fouvryten} for some discussion of earlier Bombieri--Vinogradov type results for smooth numbers. In particular, Granville~\cite{granville} proved such a result, with a bound $\ll_{A} \Psi(x,y)/\log^{A}x$, that is valid for any $100 \leq y \leq x$, but only when $Q \leq \min\{\sqrt{x}/\log^{B(A)}x , y^{C(A)\log\log y /\log\log\log y } \}$. When $y$ is small this range for $Q$ is much smaller than permitted in Theorem 1. The author is not aware of any previous results like Theorem 2 in the literature, although the methods of Fouvry and Tenenbaum~\cite{fouvryten} could presumably be adapted to yield a similar result on the range $\exp\{\log^{2/3+\epsilon}x\} \leq y \leq x$.


\vspace{12pt}
Now let us say something about the proofs of Theorems 1 and 2. As usual, on introducing Dirichlet characters the left hand side in Theorem 1 is seen to be
$$ \leq \sum_{q \leq Q} \frac{1}{\phi(q)} \sum_{\substack{\chi \; (\textrm{mod } q), \\ \chi \neq \chi_{0}}} |\Psi(x,y;\chi)|, \;\;\; \textrm{ where } \;\;\; \Psi(x,y;\chi) := \sum_{n \leq x} \chi(n)\textbf{1}_{\{n \in \mathcal{S}(y)\}}. $$
We will divide the double sum into three parts, according as the conductor $\textrm{cond}(\chi)$ satisfies:
\begin{enumerate}
\item $\textrm{cond}(\chi) \leq \min\{y^{\eta},e^{\eta \sqrt{\log x}}\}$; or

\item $\min\{y^{\eta},e^{\eta \sqrt{\log x}}\} < \textrm{cond}(\chi) \leq x^{\eta}$; or

\item $\textrm{cond}(\chi) > x^{\eta}$.
\end{enumerate}
Here $\eta > 0$ will be a certain sufficiently small constant.

In $\S 3$ and in the appendix we will prove a result that gives bounds on $|\Psi(x,y;\chi)|$ in terms of the zero-free region of the $L$-series $L(s,\chi)$. Combining this result with the classical zero-free region for Dirichlet $L$-functions will yield the following character sum estimate, which improves on Th\'{e}or\`{e}me 4 of Fouvry and Tenenbaum~\cite{fouvryten0}:
\begin{thm3}
There exist a small absolute constant $b > 0$, and a large absolute constant $M > 0$, such that the following is true. If $\log^{M}x \leq y \leq x$ is large; and $\chi$ is a non-principal Dirichlet character with conductor $r := \textup{cond}(\chi) \leq y^{b}$ and to a modulus $\leq x$; and the largest real zero $\beta=\beta_{\chi}$ of $L(s,\chi)$ is $\leq 1 - M/\log y$; then
$$ |\Psi(x,y;\chi)| \ll \Psi(x,y) \sqrt{\log x \log y} (e^{-(b \log x) \min\{1/\log r, 1-\beta \}}\log\log x + e^{-b\sqrt{\log x}} + y^{-b}) . $$
\end{thm3}
The part of the double sum on the range (i) can be bounded using Theorem 3, except for any characters for which $L(s,\chi)$ has a real zero too close to 1 (i.e. a Siegel zero, in a somewhat strong sense). However, the contribution from any such characters can also be bounded successfully, as we will show in $\S 3$.

On the range (ii), our bounds for $|\Psi(x,y;\chi)|$ will be satisfactory provided we have a fairly big zero-free region for $L(s,\chi)$. Although such a zero-free region isn't known for individual $L$-functions, zero-density estimates imply that it fails for a few at most. Since the summands in (ii) are accompanied by a factor $1/\phi(q)$, with $q > \min\{y^{\eta},e^{\eta \sqrt{\log x}}\}$, the contribution from such rogue $L$-functions will trivially be small enough.

On the range (iii) we will apply the multiplicative large sieve directly to bound our sums. It is now well understood that such a procedure will succeed provided we can decompose $|\Psi(x,y;\chi)|$ into (a small number of non-trivial) double character sums, and provided we look at characters whose conductor is large enough relative to the density of the $y$-smooth numbers less than $x$. We supply the relevant argument in $\S 4$.

The left hand side in Theorem 2 is equal to $\sum_{q \leq Q} (1/\phi(q)) \sum_{\chi \; (\textrm{mod } q), \; \chi \neq \chi_{0}} |\Psi(x,y;\chi)|^{2}$, and the arguments used to prove Theorem 1 will apply directly to bound this as well (considerably more easily when it comes to using the large sieve on the range (iii)).

\vspace{12pt}
The reader might wonder why it is necessary to use both zero-density methods, and the large sieve, to handle characters of large conductor here, when the classical Bombieri--Vinogradov theorem can be proved using either method on its own. (See e.g. Bombieri's paper~\cite{bomb} for a proof using zero-density estimates, and the third edition of Davenport's book~\cite{davenport} for a proof using the large sieve directly.) The answer is that the $y$-smooth numbers are much sparser than the primes, unless $y$ is very large, and are related to the $L$-series $L(s,\chi)$ in a more indirect, ``exponentiated'' way.

Fouvry and Tenenbaum~\cite{fouvryten} used the large sieve to prove a Bombieri--Vinogradov theorem for smooth numbers, and as we remarked this only works when $\exp\{\log^{2/3+\epsilon}x\} \leq y \leq x$. For smaller $y$ the $y$-smooth numbers are sufficiently sparse that the bound supplied by the large sieve, which is insensitive to this sparseness, becomes poor.

In the classical Bombieri--Vinogradov theorem, after applying the explicit formula for character sums over primes one is left to bound a double sum over characters and over zeros of $L(s,\chi)$. In contrast, in proving Theorems 1 and 2 the sum over zeros and the sum over characters become separated ``by an exponentiation''. The reader is referred to $\S 3$ for a proper explanation of this, but it essentially means that one needs good bounds for the {\em number of characters} $\chi$ (with certain conductors) such that $L(s,\chi)$ has at least one zero in a box, as opposed to bounds for the total number of zeros of all such $L(s,\chi)$ in that box. Zero-density estimates of the latter type imply bounds of the former type, but with a quantitative loss that is considerable when counting zeros far into the critical strip or of large height. This loss means that we cannot handle characters whose conductor is too big by zero-density methods.

\vspace{12pt}
Finally we point out that our approach of relating $\Psi(x,y;\chi)$ to zeros of $L(s,\chi)$ is connected to various previous work. Within the classical zero-free region, this idea appears to originate with Fouvry and Tenenbaum~\cite{fouvryten0}. Later, Konyagin and Soundararajan~\cite{konsound} and Soundararajan~\cite{sound} (see particularly $\S 3$ of that paper) investigated the consequences of a larger zero-free region, in conjunction with zero-density estimates. The author extended some of that work as an ingredient in the paper~\cite{harper3}. In $\S 3$ of the present paper we prove a result connecting $\Psi(x,y;\chi)$ to zeros of $L(s,\chi)$, that is roughly comparable to what is obtained for character sums over primes using Perron's formula. As we will discuss, this result is still capable of various refinements, but the author hopes it will be a useful general tool in future work on smooth numbers in arithmetic progressions.

\section{Background on smooth numbers and on zeros of $L$-functions}
Our arguments will require a few background results on the distribution of smooth numbers. Most importantly, we shall require the following seminal result of Hildebrand and Tenenbaum~\cite{ht}:
\begin{smooth1}[Hildebrand and Tenenbaum, 1986]
We have uniformly for $x \geq y \geq 2$,
$$ \Psi(x,y) = \frac{x^{\alpha} \zeta(\alpha,y)}{\alpha \sqrt{2\pi(1+(\log x)/y)\log x \log y}} \left(1 + O\left(\frac{1}{\log(u+1)} + \frac{1}{\log y} \right) \right), $$
where $u = (\log x)/\log y$, $\zeta(s,y) := \sum_{n: n \textrm{ is } y \textrm{ smooth}} 1/n^{s} = \prod_{p \leq y}(1-p^{-s})^{-1}$ for $\Re(s) > 0$, and $\alpha = \alpha(x,y) > 0$ is defined by
$$ \sum_{p \leq y} \frac{\log p}{p^{\alpha}-1} = \log x. $$
\end{smooth1}

Smooth Numbers Result 1 is proved using a saddle-point method, and we will sometimes refer to this expression for $\Psi(x,y)$ as the ``saddle-point expression''. Hildebrand and Tenenbaum~\cite{ht} also established a simple approximation for $\alpha(x,y)$ on the whole range $2 \leq y \leq x$. Their Lemma 2 implies, in particular, that when $\log x < y \leq x^{1/3}$ one has
$$ \alpha(x,y) = 1 - \frac{\log(u\log u)}{\log y} + O\left(\frac{1}{\log y}\right). $$

The reader might desire a more explicit estimate for $\Psi(x,y)$, so we remark that a result of Hildebrand~\cite{hildebrand} implies, in particular, that
$$ \Psi(x,y) = x\rho(u)\left(1+O\left(\frac{\log(u+1)}{\log y}\right) \right), \;\;\;\;\; e^{(\log\log x)^{2}} \leq y \leq x, $$
where the Dickman function $\rho(u)$ is a certain continuous function that satisfies $\rho(u) = e^{-(1+o(1))u\log u}$ as $u \rightarrow \infty$. Thus the $y$-smooth numbers are a very sparse set when $u=(\log x)/\log y$ is large. See Hildebrand and Tenenbaum's paper~\cite{ht} for much further discussion of the behaviour of $\Psi(x,y)$.

We will also require some crude ``local'' information about $\Psi(x,y)$, that describes roughly how this function changes when $x$ or $y$ change a little.
\begin{smooth2}[Following Hildebrand and Tenenbaum, and others]
For any large $\log x \leq y \leq x$ we have
$$ \Psi(2x,y) \ll \Psi(x,y) \;\;\; \textrm{ and } \;\;\; \Psi(x,y(1+1/\log x)) \ll \Psi(x,y) . $$
\end{smooth2}
The first bound follows immediately from Theorem 3 of Hildebrand and Tenenbaum~\cite{ht}, for example. For the second bound, we may assume that $y \leq x^{1/3}$ (since otherwise $\Psi(x,y) \gg x$), and in view of Smooth Numbers Result 1 it will suffice to show that $\alpha' := \alpha(x,y(1+1/\log x))$ satisfies $\alpha' = \alpha(x,y) + O(1/\log x)$ when $y \geq \log x$. But by definition of $\alpha'$ we have
$$ \sum_{p \leq y(1+1/\log x)} \frac{\log p}{p^{\alpha'}-1} = \log x , $$
and therefore
$$ \sum_{p \leq y} \frac{\log p}{p^{\alpha'}-1} = \log x - O(\frac{y^{1-\alpha'}\log y}{\log x}) = \log x - O(\log u) = \log(x/u^{O(1)}) . $$
By definition of the saddle-point $\alpha(\cdot,\cdot)$, this implies that $\alpha' = \alpha(x/u^{O(1)},y)$. Then the claimed estimate $\alpha' - \alpha = O(1/\log x)$ follows because, when $y > \log x$, we have $\partial\alpha(x,y)/\partial x = O(1/(x \log x \log y)) = O(1/(x \log x \log u))$, as shown in e.g. the proof of Theorem 4 of Hildebrand and Tenenbaum~\cite{ht}.

Finally, in some of our applications of Perron's formula we will need an upper bound for the quantity of $y$-smooth numbers in short intervals. The following result, which is a consequence of a ``sublinearity result'' of Hildebrand~\cite{hildebrandshort}, will be sufficient.
\begin{smooth3}[Hildebrand, 1985]
For any large $y,z$ and $x \geq \max\{y,z\}$ we have
$$ \Psi(x+z,y) - \Psi(x,y) \ll \Psi(z,y) \ll 2^{\log(x/z)/\log y} \left(\frac{z}{x}\right)^{\alpha(x,y)} \Psi(x,y) , $$
where $\alpha(x,y)$ is the saddle-point defined in Smooth Numbers Result 1.
\end{smooth3}
The first inequality here follows from Theorem 4 of Hildebrand~\cite{hildebrandshort}, and the second by iteratively applying Theorem 3 of Hildebrand and Tenenbaum~\cite{ht} with the choice $c = \min\{x/z,y\}$.

\vspace{12pt}
As the reader might expect having read the introduction, our arguments will also require various information about the zeros of Dirichlet $L$-functions $L(s,\chi)$. The following statement collects together the facts we will need.
\begin{zeros1}
There is an absolute and effective constant $\kappa > 0$ such that, for any $q, Q \geq 1$, the functions $F_{q}(s):=\prod_{\chi (\textrm{mod } q)} L(s,\chi)$ and $G_{Q}(s):=\prod_{q \leq Q} \prod_{\chi (\textrm{mod } q)}^{*} L(s,\chi)$ have the following properties (where $\prod^{*}$ denotes a product over primitive characters):
\begin{enumerate}
\item {\em (zero-free region)} $F_{q}(\sigma+it)$ has at most one zero in the region $\sigma \geq 1 - \kappa/\log(q(2+|t|))$. If such an ``exceptional'' zero exists then it is real, simple, and corresponds to a non-principal real character.

\item {\em (Page's theorem)} $G_{Q}(\sigma+it)$ has at most one zero (which, if it exists, is necessarily real, simple, and arises from a real character) in the region $\sigma \geq 1 - \kappa/\log(Q(2+|t|))$.

\item {\em (Siegel's theorem)} for any $\epsilon > 0$ there is a constant $C(\epsilon) > 0$, {\em which in general is non-effective}, such that $F_{q}(\sigma)$ has no real zeros $\sigma \geq 1 - C(\epsilon)/q^{\epsilon}$.

\item {\em (log-free zero-density estimate)} for any $\epsilon > 0$ and any $\sigma \geq 1/2$, $T \geq 1$, the function $F_{q}(s)$ has $\ll_{\epsilon} (qT)^{(12/5+\epsilon)(1-\sigma)}$ zeros $s$, counted with multiplicity, with $\Re(s) \geq \sigma$ and $|\Im(s)| \leq T$. Moreover, $G_{Q}(s)$ has $\ll_{\epsilon} (Q^{2}T)^{(12/5+\epsilon)(1-\sigma)}$ zeros in that region.
\end{enumerate}
\end{zeros1}

The zero-free region, Page's theorem and Siegel's theorem are all proved in standard textbooks on multiplicative number theory: see e.g. Chapter 11 of Montgomery and Vaughan~\cite{mv}. The log-free zero-density estimates stated above are proved in Huxley's paper~\cite{huxley} for values of $\sigma$ bounded away from 1, and in Jutila's paper~\cite{jutila} for values of $\sigma$ close to 1 (actually with a smaller exponent than the famous $12/5 + \epsilon$). The description ``log-free'' refers to the fact that there are no logarithmic factors, or other factors that do not decay with $1-\sigma$, in the estimates, which are therefore still very useful when $\sigma$ is very close to 1. We will exploit this state of affairs in e.g. Lemma 1, below.

\section{The zero-density argument}
In this section our main goal is to prove the following proposition, which gives Perron-type bounds for $\Psi(x,y;\chi)$ in terms of the zero-free region of $L(s,\chi)$.

\begin{prop1}
There exist a small absolute constant $d > 0$, and a large absolute constant $C > 0$, such that the following is true.

Let $\log^{1.1}x \leq y \leq x$ be large. Suppose that $\chi$ is a non-principal Dirichlet character with conductor $r:=\textup{cond}(\chi) \leq x^{d}$, and to modulus $q \leq x$, such that $L(s,\chi)$ has no zeros in the region
$$ \Re(s) > 1- \epsilon, \;\;\; |\Im(s)| \leq H, $$
where $C/\log y < \epsilon \leq \alpha(x,y)/2$ and $y^{0.9\epsilon} \log^{2}x \leq H \leq x^{d}$. Suppose, moreover, that {\em at least one} of the following holds:
\begin{itemize}
\item $y \geq (Hr)^{C}$;

\item $\epsilon \geq 40\log\log(qyH)/\log y$.
\end{itemize}

Then we have the bound
$$ |\Psi(x,y;\chi)| \ll \Psi(x,y)\sqrt{\log x \log y}(x^{-0.3\epsilon}\log H + \frac{1}{H^{0.02}}) . $$
\end{prop1}

We will then use Proposition 1 to deduce Theorem 3, and to handle characters on the ranges (i) and (ii) (for Theorems 1 and 2) described in the introduction.

Before we do this, let us make a few remarks. The restrictions on $\epsilon$ and $H$ in Proposition 1 may seem technical and off-putting, but they will be easy to satisfy in practice and, as the reader will see, they arise naturally in the proof. The restriction that $y \geq \log^{1.1}x$ is, to some extent, for simplicity, in particular so that
$$ \alpha(x,y) \geq \alpha(x,\log^{1.1}x) = 1-\frac{\log\log x + O(1)}{\log(\log^{1.1}x)} \geq 0.05 \gg 1, $$
and many of our arguments work when $y$ is smaller. However, when $y \leq \log x$ there is a genuine change in the nature of the $y$-smooth numbers less than $x$, in that almost all of them are products of large powers of primes. Moreover, when $y \leq \log x$ then $\Psi(x,y) = x^{o(1)}$, (see e.g. Chapter 7.1 of Montgomery and Vaughan~\cite{mv}), so one couldn't have a Bombieri--Vinogradov type theorem with a large range of summation over $q$.

Secondly, it would be more natural to compare $|\Psi(x,y;\chi)|$ with its trivial bound $\Psi_{q}(x,y)$, rather than $\Psi(x,y)$. We refrain from doing this because it would be quite complicated to formulate a single result that holds on a very wide range of $x,y$ and $q$, but again many of our arguments do supply bounds involving $\Psi_{q}(x,y)$, and a reader who wants such a result should have no difficulty in adapting our methods. Also see de la Bret\`{e}che and Tenenbaum's paper~\cite{dlbten} for many results relating $\Psi_{q}(x,y)$ and $\Psi(x,y)$.

Thirdly, the multiplier $\sqrt{\log x \log y}$ in Proposition 1 and Theorem 3, which will be problematic when $\epsilon$ is very small or $r$ is very close to $x$, can almost certainly be removed by adapting the ``majorant principle'' argument in $\S 2.3$ of the author's paper~\cite{harper3}. One can perhaps also remove the condition that either $y \geq (Hr)^{C}$ or $\epsilon \geq 40\log\log(qyH)/\log y$, by introducing smooth weights into the explicit formula arguments that prove Proposition 1 (as was done in less general settings by Soundararajan~\cite{sound} and the author~\cite{harper3}). However, since the proof of Proposition 1 is already quite involved we do not work these extensions out here.

\subsection{Proof of Proposition 1}
The proposition will be a relatively easy consequence of the following lemmas:
\begin{lem1}
There is a large absolute constant $C > 0$ such that the following is true. Suppose $\chi$ is a non-principal Dirichlet character with conductor $r$, and to modulus $q$. If $L(s,\chi)$ has no zeros in the region
$$ \Re(s) > 1- \epsilon, \;\;\; |\Im(s)| \leq H, $$
where $0 < \epsilon \leq 1/2$ and $H \geq 4$, then for any $z \geq (Hr)^{C}$, any $|t| \leq H/2$, and any $0 \leq \sigma < 1$ we have
$$ \left| \sum_{n \leq z} \frac{\Lambda(n) \chi(n)}{n^{\sigma + it}} \right| \ll \frac{z^{1-\sigma-0.9\epsilon}}{1-\sigma} + \frac{z^{1-\sigma} \log^{2}(qzH)}{(1-\sigma)H} + \log(rH) + \log^{0.9}q + \frac{1}{\epsilon} . $$
\end{lem1}

\begin{lem2}
Let the situation be as in Lemma 1, but with the condition that $z \geq (Hr)^{C}$ replaced by the condition $z \geq 2$. Then
$$ \left| \sum_{n \leq z} \frac{\Lambda(n) \chi(n)}{n^{\sigma + it}} \right| \ll \frac{z^{1-\sigma-0.95\epsilon}\log^{2}(qzH)}{1-\sigma} + \frac{z^{1-\sigma} \log^{2}(qzH)}{(1-\sigma)H} + \log(rH) + \log^{0.9}q + \frac{1}{\epsilon} . $$
\end{lem2}

We defer the proofs of Lemmas 1 and 2 to the appendix, but we remark that the proof of Lemma 1 itself uses a log-free zero-density estimate to remove various logarithmic factors, which is why it can supply a non-trivial bound even if $\epsilon$ is very small (unlike Lemma 2). Readers familiar with the proofs of Linnik's theorem on the least prime in an arithmetic progression should find this familiar. As we will soon see, the sums in Lemmas 1 and 2 will essentially appear as exponents in the proof of Proposition 1, and clearly a loss of logarithmic factors in an exponent would not yield acceptable bounds.

\vspace{12pt}
To deduce Proposition 1 we obtain bounds on $|\log L(\sigma+it,\chi;y) - \log L(\alpha+it,\chi;y)|$, where $\alpha=\alpha(x,y)$ is the saddle-point from Smooth Numbers Result 1, and
$$ L(s,\chi;y):=\prod_{p \leq y} \left(1-\frac{\chi(p)}{p^{s}}\right)^{-1} = \sum_{n \in \mathcal{S}(y)} \frac{\chi(n)}{n^{s}}, \;\;\; \Re(s) > 0 $$
is the Dirichlet series corresponding to the $y$-smooth numbers, and where $\alpha-0.8\epsilon \leq \sigma \leq \alpha$ and $|t| \leq H/2$. Indeed, remembering that we have $\epsilon \leq \alpha/2$ and $\alpha \gg 1$ in Proposition 1 (since $y \geq \log^{1.1}x$), this difference is certainly at most
\begin{eqnarray}
(\alpha-\sigma) \sup_{\sigma \leq \sigma' \leq \alpha} \left|\frac{L'(\sigma'+it,\chi;y)}{L(\sigma'+it,\chi;y)} \right| & \ll & (\alpha - \sigma)(\sup_{\alpha-0.8\epsilon \leq \sigma' \leq \alpha} \left|\sum_{n \leq y} \frac{\Lambda(n) \chi(n)}{n^{\sigma'+it}} \right| + \sum_{p \leq y} \frac{\log p}{p^{2(\alpha-0.8\epsilon)}} ) \nonumber \\
& \ll & (\alpha - \sigma)(\sup_{\alpha-0.8\epsilon \leq \sigma' \leq \alpha} \left|\sum_{n \leq y} \frac{\Lambda(n) \chi(n)}{n^{\sigma'+it}} \right| + \frac{y^{1-\alpha-0.1\epsilon}}{1-\alpha} + \frac{1}{\epsilon} ), \nonumber
\end{eqnarray}
since $2(\alpha-0.8\epsilon) \geq \alpha + 0.4\epsilon$. If $y \geq (Hr)^{C}$ then Lemma 1 implies this is all
$$ \ll (\alpha - \sigma)(\frac{y^{1-\alpha-0.1\epsilon}}{1-\alpha} + \frac{y^{1-\alpha+0.8\epsilon}\log^{2}(qyH)}{(1-\alpha)H} + \log(rH) + \log^{0.9}q + \frac{1}{\epsilon} ), $$
and since we assume in Proposition 1 that $H \geq y^{0.9\epsilon} \log^{2}x$ the second term may be omitted. If $y < (Hr)^{C}$ then the conditions of Proposition 1 will only be satisfied if $\epsilon \geq 40\log\log(qyH)/\log y$, in which case we can use Lemma 2 instead of Lemma 1, with the additional saving of $y^{-0.05\epsilon}$ in the first term there compensating for the multiplier $\log^{2}(qyH)$. Thus in any event we have
$$ |\log L(\sigma+it,\chi;y) - \log L(\alpha+it,\chi;y)| \ll (\alpha - \sigma)(\frac{y^{1-\alpha-0.1\epsilon}}{1-\alpha} + \log(rH) + \log^{0.9}q + \frac{1}{\epsilon} ). $$

Next, we know that $\alpha(x,y) = 1 - (\log(u\log u) + O(1))/\log y$, and so the above is
$$ \ll (\alpha-\sigma)(y^{-0.1\epsilon}u\log y + \log(rH) + \log^{0.9}q + \frac{1}{\epsilon}) = (\alpha-\sigma)(y^{-0.1\epsilon}\log x + \log(rH) + \log^{0.9}q + \frac{1}{\epsilon}) . $$
Since we have $\epsilon > C/\log y \geq C/\log x$, and $r,H \leq x^{d}$, and $q \leq x$ in Proposition 1, where $C$ is large and $d$ is small, we finally obtain that
$$ |\log L(\sigma+it,\chi;y) - \log L(\alpha+it,\chi;y)| \leq \frac{(\alpha-\sigma)\log x}{2} \;\;\;\;\; \textrm{if } \alpha-0.8\epsilon \leq \sigma \leq \alpha \;\; \textrm{and } |t| \leq H/2 . $$

\vspace{12pt}
Now we can prove Proposition 1 by expressing $\Psi(x,y;\chi)$ as a contour integral involving $L(\alpha+it,\chi;y)$, and shifting the line of integration. However, because the $y$-smooth numbers may be a very sparse set, and we only have useful information about $L(s,\chi;y)$ when $|t| \leq H/2$, we need to be careful about the truncation errors that arise in doing this. Using the truncated Perron formula, as in e.g. Theorems 5.2 and 5.3 of Montgomery and Vaughan~\cite{mv}, we find
\begin{eqnarray}
\Psi(x,y;\chi) & = & \frac{1}{2\pi i} \int_{\alpha-iH/2}^{\alpha+iH/2} L(s,\chi;y) \frac{x^{s}}{s} ds + O(1 + \frac{x^{\alpha}L(\alpha,\chi_{0};y)}{H} + \sum_{\substack{x/2 < n < 2x, \\ n \textrm{ is } y \textrm{ smooth}, \\ (n,q)=1}} \min\{1,\frac{x}{H|x-n|}\} ) \nonumber \\
& = & \frac{1}{2\pi i} \int_{\alpha-iH/2}^{\alpha+iH/2} L(s,\chi;y) \frac{x^{s}}{s} ds + O(1 + \frac{x^{\alpha}L(\alpha,\chi_{0};y)}{\sqrt{H}} + \sum_{\substack{|n-x| \leq x/\sqrt{H}, \\ n \textrm{ is } y \textrm{ smooth}, \\ (n,q)=1}} 1 ), \nonumber
\end{eqnarray}
where the second equality uses the Rankin-type upper bound
$$ \sum_{\substack{x/2 < n < 2x, \\ n \textrm{ is } y \textrm{ smooth}, \\ (n,q)=1}} 1 \ll \sum_{n : n \textrm{ is } y \textrm{ smooth}, (n,q)=1} \frac{x^{\alpha}}{n^{\alpha}} = x^{\alpha} L(\alpha,\chi_{0};y). $$
See e.g. Fouvry and Tenenbaum's paper~\cite{fouvryten0} for some exactly similar calculations. Our assumption that $y \geq \log^{1.1}x$ implies that $\alpha(x,y) = 1 - (\log(u\log u) + O(1))/\log y \geq 0.05$, say, and so Smooth Numbers Result 3 reveals that
$$ \sum_{\substack{|n-x| \leq x/\sqrt{H}, \\ n \textrm{ is } y \textrm{ smooth}, \\ (n,q)=1}} 1 \ll \Psi(x/\sqrt{H},y) \ll \Psi(x,y) (\sqrt{H})^{-0.04} = \Psi(x,y) H^{-0.02}. $$
Moreover, since we clearly have $\Psi(x/\sqrt{H},y) \gg 1$ (since we assume that $H \leq x^{d}$) this term includes the $O(1)$ term in our preceding expression for $\Psi(x,y;\chi)$. In addition, Smooth Numbers Result 1 implies $x^{\alpha}L(\alpha,\chi_{0};y) \leq x^{\alpha} \zeta(\alpha,y) \ll \Psi(x,y) \sqrt{\log x \log y}$, and so
$$ \Psi(x,y;\chi) = \frac{1}{2\pi i} \int_{\alpha-iH/2}^{\alpha+iH/2} L(s,\chi;y) \frac{x^{s}}{s} ds + O(\Psi(x,y)(\frac{\sqrt{\log x \log y}}{\sqrt{H}} + \frac{1}{H^{0.02}})) . $$

Finally, if we shift the line of integration and use our bound on $|\log L(\sigma+it,\chi;y) - \log L(\alpha+it,\chi;y)|$ we see
\begin{eqnarray}
\frac{1}{2\pi i} \int_{\alpha-iH/2}^{\alpha+iH/2} L(s,\chi;y) \frac{x^{s}}{s} ds & = & \frac{1}{2\pi i} \int_{\alpha-0.8\epsilon-iH/2}^{\alpha-0.8\epsilon+iH/2} L(s,\chi;y) \frac{x^{s}}{s} ds + \nonumber \\
&& + O(\frac{L(\alpha,\chi_{0};y)}{H} \int_{\alpha-0.8\epsilon}^{\alpha} e^{(\alpha-\sigma)(\log x)/2} x^{\sigma} d\sigma) \nonumber \\
& = & O(x^{\alpha}L(\alpha,\chi_{0};y)(x^{-0.8\epsilon+\epsilon/2}(\frac{1}{\alpha} + \log H) + \frac{1}{H})) \nonumber \\
& = & O(\Psi(x,y)\sqrt{\log x \log y}(x^{-0.3\epsilon}\log H + \frac{1}{H})) , \nonumber
\end{eqnarray}
from which the bound claimed in Proposition 1 immediately follows.
\begin{flushright}
Q.E.D.
\end{flushright}

\subsection{Proof of Theorem 3}
We will apply Proposition 1 with a suitable choice of $\epsilon$ and $H$. We assume, as we may, that the value of $b > 0$ in Theorem 3 was set small enough in terms of the values $d,C$ in Proposition 1 and the constant $\kappa$ in Zeros Result 1, and also that the value of $M$ in Theorem 3 was set large enough in terms of $b,d,C$.

In Theorem 3 we have $r := \textrm{cond}(\chi) \leq y^{b}$, so for any $2 \leq H \leq y^{50b}$ we have
$$ (Hr)^{C} \leq y^{51bC} \leq y . $$
Moreover, Zeros Result 1 (applied to the primitive character inducing $\chi$) and the assumptions of Theorem 3 imply that $L(s,\chi)$ has no zeros in the region
$$ \Re(s) > \max\{\beta_{\chi},1 - \frac{\kappa}{\log(r(2+H))} \}, \;\;\; |\Im(s)| \leq H , $$
where $\max\{\beta_{\chi},1 - \kappa/\log(r(2+H))\} \leq \max\{1-M/\log y,1 - \kappa/\log(y^{b}(2+y^{50b}))\} \leq 1 - C/\log y$. So if we choose
$$ H = \min\{y^{50b},e^{\sqrt{\log x}},e^{(\log x)/\log r}\log^{2}x , x^{1-\beta_{\chi}}\log^{2}x\}, \;\;\; \textrm{and} \;\;\; \epsilon = 1 - \max\{\beta_{\chi},1 - \frac{\kappa}{\log(r(2+H))} \}, $$
the reader can easily check that $H \geq y^{0.9\epsilon}\log^{2}x$ (bearing in mind that $y \geq \log^{M}x$ in Theorem 3), and so all the conditions of Proposition 1 will be satisfied. Since we have
$$ \epsilon \geq \min\{1-\beta_{\chi},\frac{\kappa}{2\log r},\frac{\kappa}{2\log(2+H)}\} \gg \min\{1-\beta_{\chi},\frac{1}{\log r},\frac{1}{\sqrt{\log x}}\}, $$
we conclude that
\begin{eqnarray}
|\Psi(x,y;\chi)| & \ll & \Psi(x,y)\sqrt{\log x \log y}(x^{-0.3\epsilon}\log H + \frac{1}{H^{0.02}}) \nonumber \\
& \ll & \Psi(x,y)\sqrt{\log x \log y}(x^{-0.3\epsilon}\log H + y^{-b} + e^{-b\sqrt{\log x}} + e^{-b(\log x)/\log r} + e^{-(b\log x)(1-\beta)}) \nonumber \\
& \ll & \Psi(x,y) \sqrt{\log x \log y} (e^{-(b \log x) \min\{1/\log r, 1-\beta \}}\log\log x + y^{-b} + e^{-b\sqrt{\log x}}), \nonumber
\end{eqnarray}
as asserted in Theorem 3.
\begin{flushright}
Q.E.D.
\end{flushright}

\subsection{Application to Theorems 1 and 2}
Recall that in Theorem 1 we are trying to bound
$$ \sum_{q \leq Q} \max_{(a,q)=1} \left| \Psi(x,y;q,a)-\frac{\Psi_{q}(x,y)}{\phi(q)} \right| \leq \sum_{q \leq Q} \frac{1}{\phi(q)} \sum_{\substack{\chi \; (\textrm{mod } q), \\ \chi \neq \chi_{0}}} |\Psi(x,y;\chi)|. $$
In this subsection we will show that, provided $\eta$ is fixed small enough in terms of the various constants $b,M,d,C$ in Theorem 3 and Proposition 1, and provided the constants $c,K$ in Theorem 1 are fixed suitably small and large (respectively) in terms of $\eta$, then
\begin{eqnarray}\label{zdtarget}
\sum_{1 < r \leq x^{\eta}} \sum_{\substack{\chi^{*} (\textrm{mod } r), \\ \chi^{*} \; \textrm{primitive}}} \sum_{q \leq Q} \frac{1}{\phi(q)} \sum_{\substack{\chi (\textrm{mod } q), \\ \chi^{*} \; \textrm{induces } \chi}} |\Psi(x,y;\chi)| \ll \Psi(x,y)\left(e^{-\frac{cu}{\log^{2}(u+1)}} + y^{-c}\right) ,
\end{eqnarray}
and also that it satisfies the stronger ineffective bound claimed in Theorem 1. We will also prove corresponding statements for Theorem 2. The remaining characters, with conductor $x^{\eta} < r \leq Q$, will be dealt with in $\S 4$.

We remark that we will break the sum in (\ref{zdtarget}) into a few different pieces, as suggested in the introduction, and for most of the pieces will obtain bounds of the form
$$ \ll \Psi(x,y)(e^{-\Theta(\sqrt{\log x})} + y^{-\Theta(1)}). $$
This is certainly acceptable for (\ref{zdtarget}), since (recalling that $u := (\log x)/\log y$) we always have $e^{-\sqrt{\log x}} \leq \max\{y^{-1},e^{-u}\}$. As in many multiplicative problems, a term $e^{-\sqrt{\log x}}$ arises from balancing contributions of the form $e^{-(\log x)/\log r} + r^{-1}$ and $e^{-(\log x)/\log H} + H^{-1}$ (as we have already seen in the proof of Theorem 3). The only instance in which we actually obtain a bound $ \ll \Psi(x,y)\left(e^{-cu/\log^{2}(u+1)} + y^{-c}\right)$ is, as we shall see, when considering the contribution from an exceptional character that gives rise to a Siegel zero.

For ease of writing, let us introduce some temporary notation: we will let
$$ \mathcal{G}_{1} := \bigcup_{1 < r \leq \min\{y^{\eta},e^{\eta \sqrt{\log x}}\}} \{\chi^{*} \; (\textrm{mod } r) : L(s,\chi^{*}) \; \textrm{has no real zero that is } > 1 - \frac{M}{\min\{\log y ,\sqrt{\log x}\}}\}, $$
$$ \mathcal{G}_{2} := \bigcup_{\min\{y^{\eta},e^{\eta \sqrt{\log x}}\} < r \leq x^{\eta}} \{\chi^{*} \; (\textrm{mod } r) : L(s,\chi^{*}) \neq 0 \; \textrm{for any } \Re(s) > \frac{299}{300}, \; |\Im(s)| \leq r^{100}\}, $$
where $M$ is the absolute constant from the statement of Theorem 3. These are sets of ``good'' characters corresponding to the ranges (i) and (ii) described in the introduction. Indeed, using Theorem 3 we see the contribution to (\ref{zdtarget}) from characters induced from $\chi^{*} \in \mathcal{G}_{1}$ has order at most
\begin{eqnarray}
&& \sum_{\chi^{*} \in \mathcal{G}_{1}} \sum_{q \leq Q} \frac{1}{\phi(q)} \sum_{\substack{\chi \; (\textrm{mod } q), \\ \chi^{*} \; \textrm{induces } \chi}} \Psi(x,y)\sqrt{\log x \log y}(e^{-b\sqrt{\log x}} + y^{-b}) \nonumber \\
& \ll & \Psi(x,y)\sqrt{\log x \log y}(e^{-b\sqrt{\log x}} + y^{-b}) \sum_{r \leq \min\{y^{\eta},e^{\eta \sqrt{\log x}}\}} \sum_{\substack{\chi^{*} \; (\textrm{mod } r), \\ \chi^{*} \in \mathcal{G}_{1}}} \frac{1}{\phi(r)} \sum_{s \leq Q/r} \frac{1}{\phi(s)} \nonumber \\
& \ll & \Psi(x,y)(e^{-b\sqrt{\log x}/2} + y^{-b/2}), \nonumber
\end{eqnarray}
where the final line uses the fact that $\eta$ is small in terms of $b$, and also the fact that we have $y \geq \log^{K}x$ in Theorem 1 (so we can absorb all the logarithmic factors). This is certainly an acceptable bound. Using Proposition 1 with the choices
$$ \epsilon = \min\{1/300,(10\log r)/\log y\} \;\;\; \textrm{and} \;\;\; H = r^{100}, $$
(which do satisfy the conditions $40\log\log(qyH)/\log y \leq \epsilon \leq \alpha(x,y)/2$ and $y^{0.9\epsilon}\log^{2}x \leq H \leq x^{d}$, since we have $\min\{y^{\eta},e^{\eta \sqrt{\log x}}\} < r \leq x^{\eta}$ and $y \geq \log^{K}x$ with $K$ large), the contribution to (\ref{zdtarget}) from characters induced from $\chi^{*} \in \mathcal{G}_{2}$ is also seen to be
\begin{eqnarray}
& \ll & \sum_{\min\{y^{\eta},e^{\eta \sqrt{\log x}}\} < r \leq \min\{y^{1/3000},x^{\eta}\}} \sum_{\substack{\chi^{*} \; (\textrm{mod } r), \\ \chi^{*} \; \textrm{primitive}}} \sum_{q \leq Q} \frac{1}{\phi(q)} \sum_{\substack{\chi \; (\textrm{mod } q), \\ \chi^{*} \; \textrm{induces } \chi}} \frac{\Psi(x,y) \sqrt{\log x \log y}}{r^{2}} \nonumber \\
&& + \sum_{\min\{y^{1/3000},x^{\eta}\} < r \leq x^{\eta}} \sum_{\substack{\chi^{*} \; (\textrm{mod } r), \\ \chi^{*} \; \textrm{primitive}}} \sum_{q \leq Q} \frac{1}{\phi(q)} \sum_{\substack{\chi \; (\textrm{mod } q), \\ \chi^{*} \; \textrm{induces } \chi}} \Psi(x,y) \sqrt{\log x \log y} \left(\frac{\log r}{x^{0.001}} + \frac{1}{r^{2}} \right) \nonumber \\
& \ll & \Psi(x,y) \sqrt{\log x \log y} \sum_{\min\{y^{\eta},e^{\eta \sqrt{\log x}}\} < r \leq x^{\eta}} \frac{1}{\phi(r) r^{2}} \sum_{\substack{\chi^{*} \; (\textrm{mod } r), \\ \chi^{*} \; \textrm{primitive}}} \sum_{s \leq Q/r} \frac{1}{\phi(s)} \nonumber \\
& \ll & \Psi(x,y) \log^{2}x (e^{-\eta \sqrt{\log x}} + y^{-\eta}), \nonumber
\end{eqnarray}
provided $\eta$ was set sufficiently small that $1/r^{2} \geq (\log r)/x^{0.001}$ in the above. This will be an acceptable bound, since we have $y \geq \log^{K}x$ and so the $\log^{2}x$ multiplier can be absorbed into the other terms.

Next, using the log-free zero-density estimate in Zeros Result 1 we see that, for any $R \geq 3$,
\begin{eqnarray}
\sum_{R < r \leq 2R} \sum_{\substack{\chi^{*} \; (\textrm{mod } r), \\ L(s,\chi^{*}) = 0 \; \textrm{for some} \\ \Re(s) > 299/300, \; |\Im(s)| \leq r^{100}}} \frac{1}{\phi(r)} & \ll & \frac{\log\log R}{R} \sum_{R < r \leq 2R} \sum_{\substack{\chi^{*} \; (\textrm{mod } r), \\ L(s,\chi^{*}) = 0 \; \textrm{for some} \\ \Re(s) > 299/300, \; |\Im(s)| \leq (2R)^{100}}} 1 \nonumber \\
& \ll & \frac{\log\log R}{R} (R^{102})^{(5/2)(1-(299/300))}, \nonumber
\end{eqnarray}
which is $\ll R^{-1/10}$, say. On splitting into dyadic intervals, it follows that
$$ \sum_{\min\{y^{\eta},e^{\eta \sqrt{\log x}}\} < r \leq x^{\eta}} \sum_{\substack{\chi^{*} \; (\textrm{mod } r), \\ \chi^{*} \notin \mathcal{G}_{2}}} \sum_{q \leq Q} \frac{1}{\phi(q)} \sum_{\substack{\chi \; (\textrm{mod } q), \\ \chi^{*} \; \textrm{induces } \chi}} |\Psi(x,y;\chi)| \ll \Psi(x,y) \log x (y^{-\eta/10} + e^{-\eta \sqrt{\log x}/10}), $$
which is acceptable for (\ref{zdtarget}).

At this point the only contribution to (\ref{zdtarget}) we have not dealt with is that of characters $\chi$ with conductor at most $\min\{y^{\eta},e^{\eta \sqrt{\log x}}\}$ and a real zero that is $> 1 - M/\min\{\log y ,\sqrt{\log x}\}$. Using Page's Theorem (from Zeros Result 1), provided $\eta$ was chosen small enough in terms of $M$ there will be at most one primitive character $\chi^{*}_{\textrm{bad}}$ giving rise to such contributions. If such $\chi^{*}_{\textrm{bad}}$ exists, and has conductor $r_{\textrm{bad}}$, then the contribution is
$$ \sum_{\substack{q \leq Q, \\ r_{\textrm{bad}} \mid q}} \frac{1}{\phi(q)} \sum_{\substack{\chi \; (\textrm{mod } q), \\ \chi^{*}_{\textrm{bad}} \; \textrm{induces } \chi}} |\Psi(x,y;\chi)| \ll \frac{\max_{\chi} |\Psi(x,y;\chi)|}{\phi(r_{\textrm{bad}})} \sum_{s \leq Q/r_{\textrm{bad}}} \frac{1}{\phi(s)} \ll \frac{\log x }{\phi(r_{\textrm{bad}})} \max_{\chi} |\Psi(x,y;\chi)|, $$
where the maxima are over characters $\chi$ to modulus $\leq Q$ and induced by $\chi^{*}_{\textrm{bad}}$.

Let us assume, at first, that $\log^{K}x \leq y \leq x^{1/(\log\log x)^{2}}$, say. In this case we have $u = (\log x)/\log y \geq (\log\log x)^{2}$, so it will suffice for Theorem 1 if we can show that
$$ |\Psi(x,y;\chi)| \ll \Psi(x,y) \log^{2}x \left(e^{-\frac{cu}{\log^{2}(u+1)}} + y^{-c}\right) $$
for all characters $\chi$ to modulus $\leq Q$ and induced by $\chi^{*}_{\textrm{bad}}$. The point here is that we don't need to worry about the factor $\log x /\phi(r_{\textrm{bad}})$ in the previous display, or the multiplier $\log^{2}x$, since these can be absorbed into our bound when $u$ is this large (at the cost of slightly adjusting the value of $c$).

Using the truncated Perron formula exactly as in the proof of Proposition 1,
\begin{eqnarray}
\Psi(x,y;\chi) & = & \frac{1}{2\pi i} \int_{\alpha - i/(2\log y)}^{\alpha + i/(2\log y)} L(s,\chi;y) \frac{x^{s}}{s} ds + \frac{1}{2\pi i} \int_{\alpha - iy/2}^{\alpha - i/(2\log y)} L(s,\chi;y) \frac{x^{s}}{s} ds + \nonumber \\
&& + \frac{1}{2\pi i} \int_{\alpha + i/(2\log y)}^{\alpha + iy/2} L(s,\chi;y) \frac{x^{s}}{s} ds + O\left(\frac{\Psi(x,y) \sqrt{\log x \log y}}{y^{0.02}} \right). \nonumber
\end{eqnarray}
The ``big Oh'' term is acceptably small, bearing in mind that $c \leq 0.02$ is small.

To bound the integrals we employ an argument that was used to great effect\footnote{Lemma 5.2 of Soundararajan~\cite{sound}, which we shall ultimately use, has a relatively straightforward proof using the Cauchy--Schwarz inequality and estimates for the Riemann zeta function, but it will nevertheless produce good bounds. A direct argument using information about $L(\alpha+it,\chi)$ would seem (when $|t|$ is large) to be far more complicated.} by Soundararajan~\cite{sound}. Thus $|L(\alpha+it,\chi;y)/L(\alpha,\chi_{0};y)|$ is equal to
\begin{eqnarray}
\prod_{p \leq y} \left|\frac{1-\chi(p)/p^{\alpha+it}}{1-\chi_{0}(p)/p^{\alpha}}\right|^{-1} \leq \prod_{p \leq y} \left|1 + \frac{\Re(\chi_{0}(p)-\chi(p)/p^{it})}{p^{\alpha}-\chi_{0}(p)}\right|^{-1} & = & \prod_{p \leq y, p \nmid q} \left|1 + \sum_{k=1}^{\infty} \frac{1-\Re(\chi(p)/p^{it})}{p^{k\alpha}}\right|^{-1} \nonumber \\
& \leq & e^{-\sum_{p \leq y, p \nmid q} \frac{1-\Re(\chi(p)/p^{it})}{p^{\alpha}} }, \nonumber
\end{eqnarray}
where the final inequality uses the series expansion of $\exp\{(1-\Re(\chi(p)/p^{it}))/p^{\alpha}\}$, and the fact that $0 \leq 1-\Re(\chi(p)/p^{it}) \leq 2$. Using Page's theorem, if $\chi^{*}_{\textrm{bad}}$ exists then it, and all the characters $\chi$ it induces, have order two (i.e. are real and non-principal). Therefore Lemma 5.2 of Soundararajan~\cite{sound} implies that, if $1/(2\log y) \leq |t| \leq y/2$,
$$ \sum_{p \leq y, p \nmid q} \frac{1-\Re(\chi(p)/p^{it})}{p^{\alpha}} \gg \frac{u}{\log^{2}(u+1)}. $$
When $|t| < 1/(2\log y)$, Soundararajan's lemma instead yields (keeping in mind that $\chi$ and $\chi^{*}_{\textrm{bad}}$ are real characters) that
\begin{eqnarray}
\sum_{p \leq y, p \nmid q} \frac{1-\Re(\chi(p)/p^{it})}{p^{\alpha}} \gg \sum_{p \leq y, p \nmid q} \frac{1-\chi(p)}{p^{\alpha}} & = & \sum_{p \leq y, p \nmid q} \frac{1-\chi^{*}_{\textrm{bad}}(p)}{p^{\alpha}} \nonumber \\
& \geq & \frac{1}{\log y} \left(\sum_{\sqrt{y} \leq p \leq y, p \nmid q} \frac{\log p}{p^{\alpha}} - \sum_{\sqrt{y} \leq p \leq y, p \nmid q} \frac{\chi^{*}_{\textrm{bad}}(p)\log p}{p^{\alpha}} \right) . \nonumber
\end{eqnarray}
Since we assume that $r_{\textrm{bad}} = \textrm{cond}(\chi^{*}_{\textrm{bad}})$ is at most $y^{\eta}$, with $\eta$ small, partial summation from standard estimates for $\sum_{n \leq z} \Lambda(n) \chi(n)$ (as in e.g. Theorem 11.16 and Exercise 11.3.1.2 of Montgomery and Vaughan~\cite{mv}) implies
$$ \sum_{\sqrt{y} \leq p \leq y} \frac{\log p}{p^{\alpha}} - \sum_{\sqrt{y} \leq p \leq y} \frac{\chi^{*}_{\textrm{bad}}(p)\log p}{p^{\alpha}} \gg \sum_{\sqrt{y} \leq p \leq y} \frac{\log p}{p^{\alpha}} \gg \frac{y^{1-\alpha}}{1-\alpha} \gg u\log y = \log x , $$
recalling from $\S 2$ that $\alpha(x,y) = 1 - (\log(u\log u) + O(1))/\log y$. The essential point in these calculations is that the contribution to $\sum_{n \leq z} \Lambda(n) \chi(n)$ from an exceptional real zero comes with a negative sign, so when subtracted makes a positive (i.e. a helpful) contribution to our lower bounds. We need to be careful about the contribution from primes $p$ that divide $q$, but in fact this is $\ll (\log q)/y^{\alpha/2} \ll (\log q)/\log x \ll 1$ (bearing in mind that $\alpha(x,y) \geq \alpha(x,\log^{K}x) \gg 1$), which is negligible.

In summary we have shown, as we wanted, that when $\log^{K}x \leq y \leq x^{1/(\log\log x)^{2}}$,
\begin{eqnarray}
\left|\frac{1}{2\pi i} \int_{\alpha - iy/2}^{\alpha + iy/2} L(s,\chi;y) \frac{x^{s}}{s} ds \right| & \ll & L(\alpha,\chi_{0};y) e^{-\Theta(u/\log^{2}(u+1))} x^{\alpha} \log y \nonumber \\
& \ll & e^{-\Theta(u/\log^{2}(u+1))} \Psi(x,y) \sqrt{\log x} \log^{3/2}y , \nonumber
\end{eqnarray}
where the final inequality uses Smooth Numbers Result 1.

We must still bound the contribution from $\chi^{*}_{\textrm{bad}}$ in the case where $x^{1/(\log\log x)^{2}} < y \leq x$, which is potentially difficult because of various logarithmic multipliers that occur. Fortunately we can deploy existing results of Fouvry and Tenenbaum~\cite{fouvryten0,fouvryten}, who carefully investigated the influence of exceptional zeros on $\Psi(x,y;\chi)$ when $y$ isn't too small. Indeed, since we assume that $r_{\textrm{bad}} \leq e^{\eta \sqrt{\log x}} \leq y^{\eta/\log\log x}$, say, Lemme 2.2 of Fouvry and Tenenbaum~\cite{fouvryten} implies that
\begin{eqnarray}
|\Psi(x,y;\chi)| = \left|\sum_{\substack{n \leq x, n \in \mathcal{S}(y), \\ (n,q/r_{\textrm{bad}})=1}} \chi^{*}_{\textrm{bad}}(n)\right| & \ll & \left(\sum_{d \mid (q/r_{\textrm{bad}})} \frac{1}{d^{\alpha(x,y)}}\right) \Psi(x,y) \frac{\log(r_{\textrm{bad}})}{\log y} e^{-\Theta(u/\log^{2}(u+1))} \nonumber \\
& \ll & \left(\sum_{d \mid (q/r_{\textrm{bad}})} \frac{1}{\sqrt{d}}\right) \Psi(x,y) \frac{\log(r_{\textrm{bad}})}{\log x} e^{-\Theta(u/\log^{2}(u+1))} . \nonumber
\end{eqnarray}
Here the second inequality uses the facts that $\alpha(x,y) = 1 - (\log(u\log u)+O(1))/\log y \geq 1/2$ on our range of $y$, and that $\log x = u \log y$. Actually Fouvry and Tenenbaum's result would give this with the sum replaced by $\sum_{d \mid (q/r_{\textrm{bad}})} 1$, but one can check that their proof (combined with e.g. Theorem 3 of Hildebrand and Tenenbaum~\cite{ht}) gives the stronger bound that we claimed. Thus the contribution from characters induced by $\chi^{*}_{\textrm{bad}}$ is
$$ \ll \Psi(x,y) \frac{\log(r_{\textrm{bad}})}{\log x} e^{-\Theta(u/\log^{2}(u+1))} \sum_{\substack{q \leq Q, \\ r_{\textrm{bad}} \mid q}} \frac{1}{\phi(q)} \sum_{d \mid (q/r_{\textrm{bad}})} \frac{1}{\sqrt{d}} \ll \Psi(x,y) \frac{\log(r_{\textrm{bad}})}{\phi(r_{\textrm{bad}})} e^{-\Theta(u/\log^{2}(u+1))} . $$
This is obviously acceptable for the effective part of Theorem 1, as in (\ref{zdtarget}). As usual, the stronger ineffective bound $\ll_{A} \Psi(x,y) e^{-\Theta(u/\log^{2}(u+1))} \log^{-A}x$ follows from Siegel's theorem (in Zeros Result 1), which implies that if $L(s,\chi^{*}_{\textrm{bad}})$ has a real zero that is $> 1 - M/\sqrt{\log x}$ then we must have $r_{\textrm{bad}} \gg_{A} \log^{A}x$.

\vspace{12pt}
We have now completed our treatment of all characters with conductor $\leq x^{\eta}$ in Theorem 1. One can handle the contribution from such characters to Theorem 2 using exactly the same arguments, since in Theorem 2 one needs to bound
$$ \sum_{q \leq Q} \sum_{(a,q)=1} \left| \Psi(x,y;q,a)-\frac{\Psi_{q}(x,y)}{\phi(q)} \right|^{2} = \sum_{q \leq Q} \frac{1}{\phi(q)} \sum_{\substack{\chi \; (\textrm{mod } q), \\ \chi \neq \chi_{0}}} |\Psi(x,y;\chi)|^{2}, $$
so we merely insert the squares of all our bounds for $|\Psi(x,y;\chi)|$. It remains, for Theorems 1 and 2, to handle characters with conductor $> x^{\eta}$, which we shall do in the next section using the large sieve.

\section{The large sieve argument}
In this section we will use the multiplicative large sieve to complete the proofs of Theorems 1 and 2. We will apply the large sieve in the following standard form, due (apart from the values of the constants 1 and 3 multiplying $N$ and $Q^{2}$) to Gallagher~\cite{gallagher}.
\begin{multls}[Gallagher, 1967]
For any $Q \geq 1$ and any complex numbers $(a_{n})_{n=M+1}^{M+N}$, we have
$$ \sum_{q \leq Q} \frac{q}{\phi(q)} \sum_{\substack{\chi^{*} \; (\textrm{mod } q), \\ \chi^{*} \; \textrm{primitive}}} \left| \sum_{n=M+1}^{M+N} a_{n}\chi^{*}(n) \right|^{2} \leq (N+3Q^{2}) \sum_{n=M+1}^{M+N} |a_{n}|^{2}. $$
\end{multls}

Indeed, we can finish the proof of Theorem 2 almost immediately using Multiplicative Large Sieve 1. In view of the calculations in $\S 3.3$, it remains to bound the contribution from characters with conductor $ \geq x^{\eta}$, which is
\begin{eqnarray}
&& \sum_{x^{\eta} \leq r \leq Q} \sum_{\substack{\chi^{*} \; (\textrm{mod } r), \\ \chi^{*} \; \textrm{primitive}}} \sum_{\substack{x^{\eta} \leq q \leq Q, \\ r \mid q}} \frac{1}{\phi(q)} \sum_{\substack{\chi \; (\textrm{mod } q), \\ \chi^{*} \; \textrm{induces } \chi}} \left|\sum_{n \leq x, n \in \mathcal{S}(y)} \chi^{*}(n) \sum_{d \mid (n,q/r)} \mu(d) \right|^{2} \nonumber \\
& \leq & \sum_{x^{\eta} \leq r \leq Q} \sum_{\substack{\chi^{*} \; (\textrm{mod } r), \\ \chi^{*} \; \textrm{primitive}}} \sum_{\substack{x^{\eta} \leq q \leq Q, \\ r \mid q}} \frac{1}{\phi(q)} \tau(q/r) \sum_{d \mid (q/r)} |\Psi(x/d,y;\chi^{*})|^{2} \nonumber \\
& \ll & \sum_{l=0}^{[\log(Q/x^{\eta})/\log 2]} \sum_{d \leq Q/(2^{l}x^{\eta})} \frac{\log^{2}(Q/(2^{l}x^{\eta}d) + 1) \tau(d)}{\phi(d)} \sum_{2^{l}x^{\eta} \leq r \leq 2^{l+1}x^{\eta}} \frac{1}{\phi(r)} \sum_{\substack{\chi^{*} \; (\textrm{mod } r), \\ \chi^{*} \; \textrm{primitive}}} |\Psi(x/d,y;\chi^{*})|^{2} . \nonumber
\end{eqnarray}
Here the first line uses the fact that $\chi(n)=\chi^{*}(n) \textbf{1}_{(n,q/r)=1}$; on the second line we use the Cauchy--Schwarz inequality, and write $\tau(\cdot)$ for the divisor function; and the third line follows because
$$ \sum_{\substack{x^{\eta} \leq q \leq Q, \\ rd \mid q}} \frac{\phi(r)\phi(d) \tau(q/r)}{\phi(q) \tau(d)} \leq \sum_{s \leq Q/rd} \frac{\tau(s)}{\phi(s)} \leq \sum_{a \leq Q/rd} \frac{1}{\phi(a)} \sum_{b \leq Q/rda} \frac{1}{\phi(b)} \ll \log^{2}(Q/rd + 1), $$
as in e.g. Exercise 2.1.13 of Montgomery and Vaughan~\cite{mv}. Applying Multiplicative Large Sieve 1 to the inner sums, we find the above is
\begin{eqnarray}
& \ll & \sum_{l=0}^{[\log(Q/x^{\eta})/\log 2]} \sum_{d \leq Q/(2^{l}x^{\eta})} \frac{\log^{2}(Q/(2^{l}x^{\eta}d) + 1) \tau(d)}{\phi(d)} \Psi(x/d,y) \left(\frac{x/d}{2^{l}x^{\eta}} + 2^{l}x^{\eta} \right) \nonumber \\
& \ll & \sum_{l=0}^{[\log(Q/x^{\eta})/\log 2]} \log^{2}(Q/(2^{l}x^{\eta}) + 1) \Psi(x,y) \left(\frac{x}{2^{l}x^{\eta}} + 2^{l}x^{\eta}\log^{2}(Q/(2^{l}x^{\eta}) + 1) \right) \nonumber \\
& \ll & \Psi(x,y) \left(\frac{x \log^{2}Q}{x^{\eta}} + Q \right) , \nonumber
\end{eqnarray}
say. Since we assume that $y \geq \log^{K}x$ in Theorem 2, and $\Psi(x,\log^{K}x) = x^{1-1/K+o(1)}$ for any fixed $K \geq 1$ (see Corollary 7.9 of Montgomery and Vaughan~\cite{mv}, or Smooth Numbers Result 1), the first term here will be $\ll \Psi(x,y)^{2}x^{-\eta/2}$ provided $K$ was set large enough in terms of $\eta$. This bound is acceptable for Theorem 2, provided the value of $c$ there is set small enough in terms of $\eta$.
\begin{flushright}
Q.E.D.
\end{flushright}

The completion of the proof of Theorem 1 will be a bit more complicated, since we do not a priori have any squares of character sums around, and it will require some care to introduce these in a way that does not spoil the resulting bounds. In $\S\S 4.1-4.3$ we shall prove the following result:
\begin{prop2}
Let $0 < \eta \leq 1/80$ be any fixed constant. Then for any large $y \leq x^{9/10}$ and any $x^{\eta} \leq Q \leq \sqrt{x}$, we have
$$ \sum_{x^{\eta} \leq r \leq Q} \sum_{\substack{\chi^{*} \; (\textrm{mod } r), \\ \chi^{*} \; \textrm{primitive}}} \sum_{x^{\eta} \leq q \leq Q} \frac{1}{\phi(q)} \sum_{\substack{\chi \; (\textrm{mod } q), \\ \chi^{*} \; \textrm{induces } \chi}} |\Psi(x,y;\chi)| \ll \log^{7/2}x \sqrt{\Psi(x,y)} \left(Q + x^{1/2 - \eta}\log^{2}x \right). $$
\end{prop2}

The reader may easily check that, together with the calculations in $\S 3.3$ (and on setting the value of $\eta$ as in $\S 3.3$), Proposition 2 will complete the proof of Theorem 1 for all $\log^{K}x \leq y \leq x^{9/10}$, provided the values of $K$ and $c$ in Theorem 1 were set suitably in terms of $\eta$.

If $x^{9/10} < y \leq x$ then the bound claimed in Proposition 2 is still true, but the arguments needed to prove this are a bit different. Indeed, on this range of $y$ the problem essentially reduces to bounding averages of character sums over primes, as in the classical Bombieri--Vinogradov theorem, and one needs to use Vaughan's Type I/Type II sums identity before applying the large sieve. We sketch a suitable argument in $\S 4.4$, which will complete the proof\footnote{We could also refer the reader to earlier Bombieri--Vinogradov type results for smooth numbers, which certainly cover the range $x^{9/10} < y \leq x$. However, the quantitative bounds in those results are not as precise as claimed in Theorem 1, and in any case it seems desirable to give a self contained treatment.} of Theorem 1 on the whole range $\log^{K}x \leq y \leq x$.

\subsection{Factoring the $y$-smooth numbers}
We can ``reveal'' any number $n$ less than $x$ by exposing its prime factors one at a time, in non-increasing order of their size. If we know that $n$ is $y$-smooth, and if $n$ isn't too small and if $y \leq x^{37/40}$, say, then none of the factors we expose will be extremely large, and so at some point as we reveal them we will have split $n$ as a product of two fairly large factors. Such an approach has been used by many previous authors, and it will allow us to decompose character sums over $y$-smooth numbers as double character sums\footnote{The reader should be forewarned that, although we use the letter $x$ in this subsection and the next, we will ultimately apply our results with $x$ replaced by $x/d$, for various small values of $d$. This is why we only postulate here that $y \leq x^{37/40}$, which is slightly weaker than the condition $y \leq x^{9/10}$ assumed in Proposition 2.}. Then we will apply the Cauchy--Schwarz inequality and the large sieve.

More precisely, let us write $P(t), p_{1}(t)$ for the greatest and least prime factors of $t \in \N$, respectively, and note that we have
\begin{eqnarray}
\Psi(x,y;\chi) & = & \Psi(x^{1/20},y;\chi) + \sum_{\substack{x^{1/20} < n \leq x, \\ n \in \mathcal{S}(y)}} \chi(n) \nonumber \\
& = & \Psi(x^{1/20},y;\chi) + \sum_{\substack{x^{1/20} < m \leq yx^{1/20}, \\ m \in \mathcal{S}(y), m/p_{1}(m) \leq x^{1/20}}} \sum_{\substack{n \leq x/m, \\ P(n) \leq p_{1}(m)}} \chi(mn) \nonumber \\
& = & \Psi(x^{1/20},y;\chi) + \sum_{i=0}^{[\log y /\log 2]} \sum_{j=0}^{[\log y / \log\lambda]} \sum_{\substack{m \in \mathcal{S}(y), m/p_{1}(m) \leq x^{1/20}, \\ 2^{i}x^{1/20} < m \leq 2^{i+1}x^{1/20}, \\ y/\lambda^{j+1} < p_{1}(m) \leq y/\lambda^{j}}} \sum_{\substack{n \leq x/m, \\ P(n) \leq p_{1}(m)}} \chi(mn) , \nonumber
\end{eqnarray}
where we set $\lambda := 1+1/(1000\log x)$. Since we assume that $y \leq x^{37/40}$ we have $yx^{1/20} \leq x^{39/40}$, and so the sums over $n \leq x/m$ are always quite long. This will be important in various of our later calculations, and particularly in $\S 4.3$ when we apply the Cauchy--Schwarz inequality and the large sieve. 

The above essentially completes the ``factorisation'' step in the proof of Proposition 2, since we have decomposed $\Psi(x,y;\chi)$ as the sum of $\Psi(x^{1/20},y;\chi)$, which will make a negligible contribution, and of a small number of double sums over $m$ and $n$. At present there is some dependence between the ranges of summation over $m$ and $n$, but we will deal with that using Perron's formula in the next subsection. However, to simplify that process we will first modify the sums a little. For the sake of concision, let us write $\mathcal{M}_{i,j}=\mathcal{M}_{i,j,x,y}$ for the range of summation in the sum over $m$. Then the quadruple sum in the previous display can be rewritten as
$$ \sum_{i=0}^{[\log y /\log 2]} \sum_{j=0}^{[\log y /\log\lambda]} \sum_{m \in \mathcal{M}_{i,j}} \chi(m) \left(\sum_{\substack{n \leq x^{19/20}/2^{i+1}, \\ P(n) \leq y/\lambda^{j+1}}} \chi(n) + \sum_{\substack{n \leq x^{19/20}/2^{i+1}, \\ y/\lambda^{j+1} < P(n) \leq y/\lambda^{j}}} \chi(n) \textbf{1}_{P(n) \leq p_{1}(m)} + \right. $$
$$ \left. + \sum_{\substack{x^{19/20}/2^{i+1} < n \leq x^{19/20}/2^{i}, \\ P(n) \leq y/\lambda^{j+1}}} \chi(n) \textbf{1}_{mn \leq x} + \sum_{\substack{x^{19/20}/2^{i+1} < n \leq x^{19/20}/2^{i}, \\ y/\lambda^{j+1} < P(n) \leq y/\lambda^{j}}} \chi(n) \textbf{1}_{mn \leq x} \textbf{1}_{P(n) \leq p_{1}(m)} \right), $$
where $\textbf{1}$ denotes the indicator function. The point of this additional decomposition is that whenever an indicator function $\textbf{1}_{V \leq W}$ appears, with $V=mn$ and $W=x$ or with $V=P(n)$ and $W=p_{1}(m)$, the quantities $V$ and $W$ are already forced to be of comparable magnitude. This will be useful in the next subsection. For the sake of concision, at some points we will use $\mathcal{N}_{i,j}^{(1)},...,\mathcal{N}_{i,j}^{(4)}$ to denote the ranges of summation over $n$ in the sums above.

We remark that the reason for dividing the values of $m$ into dyadic intervals, whilst dividing the values of $p_{1}(m)$ using the finer parameter $\lambda$, is because $\Psi(x,y)$ is more sensitive to changes in $y$ than in $x$. Recall, for example, Smooth Numbers Result 2, which we will shortly make use of.

\subsection{Separating the factors}
Using the truncated Perron formula (as in e.g. Theorems 5.2 and 5.3 of Montgomery and Vaughan~\cite{mv}), if $T > 0$ and if we set $\tilde{x} := [x] + 1/2 \in \N + 1/2$, with $[x]$ denoting the integer part of $x$, then we have
$$ \textbf{1}_{mn \leq x} = \textbf{1}_{mn < \tilde{x}} = \frac{1}{2\pi i} \int_{1/2-iT}^{1/2+iT} \frac{\tilde{x}^{s}}{m^{s}n^{s}} \frac{ds}{s} + O\left(\frac{1}{T} \left(\frac{1}{|\log(\tilde{x}/mn)|} + \sqrt{\frac{\tilde{x}}{mn}} \right) \right) . $$
Exactly similarly, we have
\begin{eqnarray}
\textbf{1}_{P(n) \leq p_{1}(m)} = \textbf{1}_{P(n) < p_{1}(m) + 1/2} & = & \frac{1}{2\pi i} \int_{1/2-iT}^{1/2+iT} \frac{(p_{1}(m)+1/2)^{s}}{P(n)^{s}} \frac{ds}{s} + \nonumber \\
&& + O\left(\frac{1}{T} \left(\frac{1}{|\log((p_{1}(m)+1/2)/P(n))|} + \sqrt{\frac{p_{1}(m)}{P(n)}} \right) \right) . \nonumber
\end{eqnarray}
In particular, if we choose $T=x^{5}$ then, since $|\tilde{x}/mn - 1| \geq 1/2mn \gg 1/x$ in the setting of $\S 4.1$, similarly for $|(p_{1}(m)+1/2)/P(n) - 1|$, both of the ``big Oh'' terms above will be $O(1/x^{4})$.

As usual, the point of applying Perron's formula in this manner is that it separates the $m$ and $n$ variables in a multiplicative way. For example, we see
\begin{eqnarray}
&& \sum_{m \in \mathcal{M}_{i,j}} \chi(m) \sum_{n \in \mathcal{N}_{i,j}^{(4)}} \chi(n) \textbf{1}_{mn \leq x} \textbf{1}_{P(n) \leq p_{1}(m)} \nonumber \\
& = & \sum_{m \in \mathcal{M}_{i,j}} \chi(m) \sum_{\substack{x^{19/20}/2^{i+1} < n \leq x^{19/20}/2^{i}, \\ y/\lambda^{j+1} < P(n) \leq y/\lambda^{j}}} \chi(n) \textbf{1}_{mn \leq x} \textbf{1}_{P(n) \leq p_{1}(m)} \nonumber \\
& = & \frac{1}{(2\pi i)^{2}} \int_{1/2-ix^{5}}^{1/2+ix^{5}} \int_{1/2-ix^{5}}^{1/2+ix^{5}} \left( \sum_{m \in \mathcal{M}_{i,j}} \chi(m) \frac{\tilde{x}^{s} (p_{1}(m)+1/2)^{u}}{m^{s}} \sum_{\substack{\frac{x^{19/20}}{2^{i+1}} < n \leq \frac{x^{19/20}}{2^{i}}, \\ y/\lambda^{j+1} < P(n) \leq y/\lambda^{j}}} \frac{\chi(n)}{n^{s}P(n)^{u}} \right) \frac{du}{u} \frac{ds}{s} \nonumber \\
&& + O(\frac{1}{x^{4}} \sum_{m \in \mathcal{M}_{i,j}} \sum_{\substack{x^{19/20}/2^{i+1} < n \leq x^{19/20}/2^{i}, \\ y/\lambda^{j+1} < P(n) \leq y/\lambda^{j}}} 1) , \nonumber
\end{eqnarray}
where crucially there is no interaction between the inner sums over $m$ and $n$.

Summarising all of the calculations from $\S\S 4.1-4.2$, we conclude that $|\Psi(x,y;\chi)|$ is
\begin{eqnarray}\label{separatedsums}
& \ll & |\Psi(x^{1/20},y;\chi)| + \sum_{i,j} \sum_{k=1}^{4} \int_{1/2-ix^{5}}^{1/2+ix^{5}} \int_{1/2-ix^{5}}^{1/2+ix^{5}} \left|\sum_{m \in \mathcal{M}_{i,j}} \chi(m) a_{s,u}^{(k)}(m) \right| \left|\sum_{n \in \mathcal{N}_{i,j}^{(k)}} \chi(n) b_{s,u}^{(k)}(n) \right| \frac{d|u|}{|u|} \frac{d|s|}{|s|} \nonumber \\
&& + \frac{\Psi(x,y)}{x^{4}} ,
\end{eqnarray}
for certain coefficients $a_{s,u}^{(k)}(m)$ and $b_{s,u}^{(k)}(n)$ whose precise forms the reader may readily ascertain, and where the ranges of summation over $i$ and $j$ are from $0$ to $[\log y / \log 2]$ and from $0$ to $[\log y / \log \lambda]$, respectively. The bound for the final ``big Oh'' term follows because, recalling that $2^{i}x^{1/20} < m \leq 2^{i+1}x^{1/20}$ and $y/\lambda^{j+1} < p_{1}(m) \leq y/\lambda^{j}$ for all $m \in \mathcal{M}_{i,j}$, we have
\begin{eqnarray}\label{smooth2app}
\sum_{i,j,k} \sum_{m \in \mathcal{M}_{i,j}} \sum_{n \in \mathcal{N}_{i,j}^{(k)}} 1 \ll \sum_{i,j} \sum_{m \in \mathcal{M}_{i,j}} \sum_{\substack{n \leq x^{19/20}/2^{i}, \\ P(n) \leq y/\lambda^{j}}} 1 \ll \sum_{i,j} \sum_{m \in \mathcal{M}_{i,j}} \sum_{\substack{n \leq x/m, \\ P(n) \leq p_{1}(m)}} 1 \leq \Psi(x,y),
\end{eqnarray}
where the second inequality uses Smooth Numbers Result 2. Actually Smooth Numbers Result 2 doesn't apply to pairs $(i,j)$ where $y/\lambda^{j} \leq \log(x^{19/20}/2^{i})$, but in that case if $y/\lambda^{j+1} < p_{1}(m) \leq y/\lambda^{j}$ then we must have $p_{1}(m) = [y/\lambda^{j}]$ anyway, since $\lambda$ is very close to 1. As we have a denominator of $x^{4}$, we could of course rely on much cruder arguments at this point, but we will require the precise calculations that we just performed in the next subsection. 

In the above, the coefficients $a_{s,u}^{(k)}(m)$ have the same order of magnitude for all $m \in \mathcal{M}_{i,j}$, for given $i,j,k$ (similarly for $b_{s,u}^{(k)}(n)$), and their products satisfy
$$ |a_{s,u}^{(k)}(m)||b_{s,u}^{(k)}(n)| \ll 1 \;\;\; \forall m \in \mathcal{M}_{i,j}, n \in \mathcal{N}_{i,j}^{(k)}, \; s, u \in [1/2-ix^{5},1/2+ix^{5}], $$
for given $i,j,k$. These properties only hold because we split our sum over $n$ into the subsums over $\mathcal{N}_{i,j}^{(k)}$ in the previous subsection, so that in our applications of Perron's formula the quantities $x$ and $mn$, and $p_{1}(m)$ and $P(n)$, are always of comparable size.

\subsection{Proof of Proposition 2}
Now we shall prove Proposition 2, which we remind the reader will complete the proof of Theorem 1 for all $\log^{K}x \leq y \leq x^{9/10}$. Exactly as at the beginning of $\S 4$, when we proved Theorem 2, we find the left hand side in Proposition 2 is
\begin{eqnarray}
&& \sum_{x^{\eta} \leq r \leq Q} \sum_{\substack{\chi^{*} \; (\textrm{mod } r), \\ \chi^{*} \; \textrm{primitive}}} \sum_{\substack{x^{\eta} \leq q \leq Q, \\ r \mid q}} \frac{1}{\phi(q)} \sum_{\substack{\chi \; (\textrm{mod } q), \\ \chi^{*} \; \textrm{induces } \chi}} \left|\sum_{n \leq x, n \in \mathcal{S}(y)} \chi^{*}(n) \sum_{d \mid (n,q/r)} \mu(d) \right| \nonumber \\
& \leq & \sum_{x^{\eta} \leq r \leq Q} \sum_{\substack{\chi^{*} \; (\textrm{mod } r), \\ \chi^{*} \; \textrm{primitive}}} \sum_{\substack{x^{\eta} \leq q \leq Q, \\ r \mid q}} \frac{1}{\phi(q)} \sum_{d \mid (q/r)} |\Psi(x/d,y;\chi^{*})| \nonumber \\
& \ll & \sum_{l=0}^{[\log(Q/x^{\eta})/\log 2]} \sum_{d \leq Q/(2^{l}x^{\eta})} \frac{\log(Q/(2^{l}x^{\eta}d) + 1)}{\phi(d)} \sum_{2^{l}x^{\eta} \leq r \leq 2^{l+1}x^{\eta}} \frac{1}{\phi(r)} \sum_{\substack{\chi^{*} \; (\textrm{mod } r), \\ \chi^{*} \; \textrm{primitive}}} |\Psi(x/d,y;\chi^{*})| . \nonumber
\end{eqnarray}

Next we would like to apply the argument of $\S\S 4.1-4.2$ to study  $\Psi(x/d,y;\chi^{*})$, but this will be problematic if $d$ is too large, because then our condition that $y \leq (x/d)^{37/40}$ may be violated. However, combining the Cauchy--Schwarz inequality with Multiplicative Large Sieve 1 yields
\begin{eqnarray}
\sum_{2^{l}x^{\eta} \leq r \leq 2^{l+1}x^{\eta}} \frac{1}{\phi(r)} \sum_{\substack{\chi^{*} \; (\textrm{mod } r), \\ \chi^{*} \; \textrm{primitive}}} |\Psi(x/d,y;\chi^{*})| & \ll & \sqrt{\sum_{2^{l}x^{\eta} \leq r \leq 2^{l+1}x^{\eta}} \frac{2^{l}x^{\eta}}{\phi(r)} \sum_{\substack{\chi^{*} \; (\textrm{mod } r), \\ \chi^{*} \; \textrm{primitive}}} |\Psi(x/d,y;\chi^{*})|^{2}} \nonumber \\
& \ll & (\sqrt{x/d} + 2^{l}x^{\eta})\sqrt{\Psi(x/d,y)}, \nonumber
\end{eqnarray}
and therefore
\begin{eqnarray}
&& \sum_{l=0}^{[\log(Q/x^{\eta})/\log 2]} \sum_{x^{1/40} < d \leq Q/(2^{l}x^{\eta})} \frac{\log(Q/(2^{l}x^{\eta}d) + 1)}{\phi(d)} \sum_{2^{l}x^{\eta} \leq r \leq 2^{l+1}x^{\eta}} \frac{1}{\phi(r)} \sum_{\substack{\chi^{*} \; (\textrm{mod } r), \\ \chi^{*} \; \textrm{primitive}}} |\Psi(x/d,y;\chi^{*})| \nonumber \\
& \ll & \sqrt{\Psi(x,y)} \sum_{l=0}^{[\log(Q/x^{\eta})/\log 2]} \log(Q/(2^{l}x^{\eta}) + 1) \sum_{x^{1/40} < d \leq Q/(2^{l}x^{\eta})} \frac{1}{\phi(d)} (\sqrt{x/d} + 2^{l}x^{\eta}) \nonumber \\
& \ll & \sqrt{\Psi(x,y)}(x^{1/2-1/80}\log^{2}Q + Q). \nonumber
\end{eqnarray}
This bound is acceptable for Proposition 2. On the other hand, when $d \leq x^{1/40}$ we have $y \leq x^{9/10} \leq (x/d)^{36/39} < (x/d)^{37/40}$, so we shall be able to apply the argument of $\S\S 4.1-4.2$.

Indeed, combining the bound (\ref{separatedsums}) with the Cauchy--Schwarz inequality and Multiplicative Large Sieve 1, we see
\begin{eqnarray}
&& \sum_{2^{l}x^{\eta} \leq r \leq 2^{l+1}x^{\eta}} \frac{1}{\phi(r)} \sum_{\substack{\chi^{*} \; (\textrm{mod } r), \\ \chi^{*} \; \textrm{primitive}}} |\Psi(x/d,y;\chi^{*})| \nonumber \\
& \ll & \frac{2^{l}x^{\eta}}{(x/d)^{3}} + \sqrt{\sum_{2^{l}x^{\eta} \leq r \leq 2^{l+1}x^{\eta}} \frac{2^{l}x^{\eta}}{\phi(r)} \sum_{\substack{\chi^{*} \; (\textrm{mod } r), \\ \chi^{*} \; \textrm{primitive}}} |\Psi((x/d)^{1/20},y;\chi^{*})|^{2}} +  \nonumber \\
&& + \frac{1}{2^{l}x^{\eta}} \sum_{i,j,k} \int \sqrt{ \sum_{r=1}^{2^{l+1}x^{\eta}} \frac{r}{\phi(r)} \sum_{\chi^{*}} \left|\sum_{m} \chi^{*}(m) a_{s,u}^{(k)}(m) \right|^{2} \cdot \sum_{r=1}^{2^{l+1}x^{\eta}} \frac{r}{\phi(r)} \sum_{\chi^{*}} \left|\sum_{n} \chi^{*}(n) b_{s,u}^{(k)}(n) \right|^{2} } \frac{d|u| d|s|}{|u||s|} \nonumber \\
& \ll & 2^{l}x^{\eta} + \sqrt{((x/d)^{1/20} + (2^{l}x^{\eta})^{2}) \Psi((x/d)^{1/20},y)} + \nonumber \\
&& + \frac{1}{2^{l}x^{\eta}} \sum_{i,j,k} \int \sqrt{(2^{i}(\frac{x}{d})^{\frac{1}{20}} + (2^{l}x^{\eta})^{2}) \sum_{m \in \mathcal{M}_{i,j}} |a_{s,u}^{(k)}(m)|^{2} (\frac{(x/d)^{\frac{19}{20}}}{2^{i}} + (2^{l}x^{\eta})^{2}) \sum_{n \in \mathcal{N}_{i,j}^{(k)}} |b_{s,u}^{(k)}(n)|^{2}} \frac{d|u| d|s|}{|u||s|}. \nonumber
\end{eqnarray}
In view of the discussion in the final paragraph of $\S 4.2$, for any $i,j,k,s,u$ we have
$\sum_{m \in \mathcal{M}_{i,j}} |a_{s,u}^{(k)}(m)|^{2} \cdot \sum_{n \in \mathcal{N}_{i,j}^{(k)}} |b_{s,u}^{(k)}(n)|^{2} \ll \#\mathcal{M}_{i,j} \#\mathcal{N}_{i,j}^{(k)}$. If we insert this upper bound then none of the terms inside the squareroot depend on $s$ or $u$ any longer, so we can perform the integrations over those variables and pick up an additional factor of $\log^{2}(x/d)$. Then applying the Cauchy--Schwarz inequality to the sum over $i,j,k$, we see the third term in the above is
\begin{eqnarray}\label{largesievepenult}
& \ll & \frac{\log^{2}(x/d)}{2^{l}x^{\eta}} \sqrt{ \sum_{i,j,k} (x/d + (2^{l}x^{\eta})^{4} + (2^{l}x^{\eta})^{2}2^{i}(x/d)^{1/20} + (2^{l}x^{\eta})^{2}\frac{(x/d)^{19/20}}{2^{i}}) \cdot \sum_{i,j,k} \#\mathcal{M}_{i,j} \#\mathcal{N}_{i,j}^{(k)} } \nonumber \\
& \ll & \log^{5/2}(x/d) \log y \left(\frac{\sqrt{x/d}}{2^{l}x^{\eta}} + 2^{l}x^{\eta} + \sqrt{y} (x/d)^{1/40} + (x/d)^{19/40} \right) \sqrt{\Psi(x/d,y)},
\end{eqnarray}
bearing in mind that the ranges of summation over $i,j,k$ are from 0 to $[\log y /\log 2]$, from 0 to $[\log y /\log \lambda] = O(\log y \log(x/d))$, and from 1 to 4 respectively. Here we also note that $\sum_{i,j,k} \#\mathcal{M}_{i,j} \#\mathcal{N}_{i,j}^{(k)} \ll \Psi(x/d,y)$, as we saw in the calculations (\ref{smooth2app}) with $x$ replaced by $x/d$.

It is clear that the bound (\ref{largesievepenult}) is also an upper bound for the two other terms in our bound for $\sum_{2^{l}x^{\eta} \leq r \leq 2^{l+1}x^{\eta}} (1/\phi(r)) \sum_{\chi^{*} (\textrm{mod } r), \; \chi^{*} \; \textrm{primitive}} |\Psi(x/d,y;\chi^{*})|$. Moreover, since we assume in Proposition 2 that $y \leq x^{9/10}$, and therefore $\sqrt{y}x^{1/40} \leq x^{19/40}$, we can replace (\ref{largesievepenult}) by the simplified upper bound
$$ \sum_{2^{l}x^{\eta} \leq r \leq 2^{l+1}x^{\eta}} \frac{1}{\phi(r)} \sum_{\substack{\chi^{*} \; (\textrm{mod } r), \\ \chi^{*} \; \textrm{primitive}}} |\Psi(x/d,y;\chi^{*})| \ll \log^{7/2}x \sqrt{\Psi(x,y)} \left(\frac{\sqrt{x/d}}{2^{l}x^{\eta}} + 2^{l}x^{\eta} + \frac{x^{19/40}}{d^{1/40}} \right) . $$
Thus the left hand side in Proposition 2 is
\begin{eqnarray}
& \ll & \log^{7/2}x \sqrt{\Psi(x,y)} \sum_{l=0}^{[\log(Q/x^{\eta})/\log 2]} \log(Q/2^{l}x^{\eta} + 1) \left(\frac{\sqrt{x}}{2^{l}x^{\eta}} + 2^{l}x^{\eta}\log(Q/2^{l}x^{\eta} + 1) + x^{19/40}\right) \nonumber \\
& \ll & \log^{7/2}x \sqrt{\Psi(x,y)} \left(x^{1/2-\eta}\log Q + Q + x^{19/40}\log^{2}Q \right) . \nonumber
\end{eqnarray}
The proposition follows immediately, given our hypothesis that $\eta \leq 1/80 < 1/40$.
\begin{flushright}
Q.E.D.
\end{flushright}

\subsection{The large sieve argument for very large $y$}
In this subsection we will sketch a proof that Proposition 2 still holds when $x^{9/10} < y \leq x$. This will finally complete the proof of Theorem 1. Arguing as in $\S 4.3$, the reader may check it will suffice to show that, for any $0 \leq l \leq \log(Q/x^{\eta})/\log 2$ and any $d \leq Q \leq \sqrt{x}$,
\begin{eqnarray}\label{largeylargesieve}
\sum_{2^{l}x^{\eta} \leq r \leq 2^{l+1}x^{\eta}} \frac{1}{\phi(r)} \sum_{\substack{\chi^{*} \; (\textrm{mod } r), \\ \chi^{*} \; \textrm{primitive}}} |\Psi(x/d,y;\chi^{*})| \ll (\frac{\sqrt{x/d}}{2^{l}x^{\eta}} + 2^{l}x^{\eta} + x^{9/20})\sqrt{x}\log^{7/2}x .
\end{eqnarray}

Firstly, if $d \geq x/y$ then $\Psi(x/d,y;\chi^{*}) = \sum_{n \leq x/d} \chi^{*}(n) = O(\sqrt{r}\log r) = O(\sqrt{x}\log x)$, in view of the P\'{o}lya--Vinogradov inequality (see e.g. Theorem 9.18 of Montgomery and Vaughan~\cite{mv}). This bound is certainly acceptable for (\ref{largeylargesieve}).

On the other hand, if $d < x/y$ then
\begin{eqnarray}
\Psi(\frac{x}{d},y;\chi^{*}) & = & \sum_{n \leq x/d} \chi^{*}(n) - \sum_{y < p \leq x/d} \sum_{m \leq x/dp} \chi^{*}(mp) \nonumber \\
& = & \sum_{n \leq x/d} \chi^{*}(n) - \sum_{m \leq x/dy} \chi^{*}(m) \sum_{y < n \leq \frac{x}{dm}} \frac{\chi^{*}(n) \Lambda(n)}{\log n} + \sum_{k=2}^{[\frac{\log x}{\log 2}]} \sum_{m \leq x/dy} \chi^{*}(m) \sum_{y < p^{k} \leq \frac{x}{dm}} \frac{\chi^{*}(p^{k})}{k} . \nonumber
\end{eqnarray}
Moreover, when $k \geq 3$ we trivially have
$$ \sum_{m \leq x/dy} \chi^{*}(m) \sum_{y < p^{k} \leq \frac{x}{dm}} \frac{\chi^{*}(p^{k})}{k} = O(\sum_{m \leq x/dy} \left(\frac{x}{dm}\right)^{1/3}) = O\left(\frac{x}{dy^{2/3}}\right) = O(\sqrt{x}), $$
and the P\'{o}lya--Vinogradov inequality again yields $\sum_{n \leq x/d} \chi^{*}(n) = O(\sqrt{x}\log x)$. Thus
\begin{eqnarray}
|\Psi(x/d,y;\chi^{*})| & \ll & \sqrt{x}\log x + \left|\sum_{m \leq x/dy} \chi^{*}(m) \sum_{y < n \leq \frac{x}{dm}} \frac{\chi^{*}(n) \Lambda(n)}{\log n}\right| + \left|\sum_{m \leq x/dy} \chi^{*}(m) \sum_{y < p^{2} \leq \frac{x}{dm}} \chi^{*}(p^{2})\right| \nonumber \\
& \ll & \sqrt{x}\log x + \frac{1}{\log y} \left|\sum_{m \leq x/dy} \chi^{*}(m) \sum_{y < n \leq \frac{x}{dm}} \chi^{*}(n) \Lambda(n) \right| + \nonumber \\
&& + \int_{y}^{x/d} \frac{1}{t\log^{2}t} \left|\sum_{m \leq x/dy} \chi^{*}(m) \sum_{t < n \leq \frac{x}{dm}} \chi^{*}(n) \Lambda(n) \right|dt + \left|\sum_{\substack{mp^{2} \leq x/d, \\ p^{2} > y}} \chi^{*}(mp^{2})\right| , \nonumber
\end{eqnarray}
where the second inequality follows by writing $1/\log n = 1/\log y - \int_{y}^{n} dt/(t\log^{2}t)$, i.e. by using partial summation.

As before, the contribution to (\ref{largeylargesieve}) from the $\sqrt{x}\log x$ term is acceptable. Bounding the last term trivially will not quite be satisfactory, but since a number less than $x$ has at most one representation as $mp^{2}$ with $p^{2} > y$, and since $\sum_{mp^{2} \leq x/d, p^{2} > y} 1 \ll x/(d\sqrt{y})$, then the Cauchy--Schwarz inequality and Multiplicative Large Sieve 1 imply that
\begin{eqnarray}
\sum_{2^{l}x^{\eta} \leq r \leq 2^{l+1}x^{\eta}} \frac{1}{\phi(r)} \sum_{\substack{\chi^{*} \; (\textrm{mod } r), \\ \chi^{*} \; \textrm{primitive}}} \left|\sum_{\substack{mp^{2} \leq x/d, \\ p^{2} > y}} \chi^{*}(mp^{2})\right| & \ll & \sqrt{\sum_{2^{l}x^{\eta} \leq r \leq 2^{l+1}x^{\eta}} \frac{2^{l}x^{\eta}}{\phi(r)} \sum_{\substack{\chi^{*} \; (\textrm{mod } r), \\ \chi^{*} \; \textrm{primitive}}} \left|\sum_{\substack{mp^{2} \leq x/d, \\ p^{2} > y}} \chi^{*}(mp^{2})\right|^{2} } \nonumber \\
& \ll & (2^{l}x^{\eta} + \sqrt{x/d}) \sqrt{x/(d\sqrt{y})}. \nonumber
\end{eqnarray}
This is acceptable for (\ref{largeylargesieve}) with much room to spare.

Finally, to bound the terms involving $\sum_{n} \chi^{*}(n)\Lambda(n)$ one can use Vaughan's identity to expand these sums into non-trivial double sums, and then collect the sum over $m \leq x/dy$ (which is a short sum, since $y$ is so close to $x$) with one of those sums and apply the Cauchy--Schwarz inequality and Multiplicative Large Sieve 1. We do not write out the details, since this is cumbersome, but refer the reader to pages 166--167 of Davenport's book~\cite{davenport} for an argument that can easily be adapted to our purposes. More specifically, one can follow that argument with the simple choices $U=V=x^{1/10}$, and discover that for the analogues of the sums $S_{1},S_{2}',S_{2}'',S_{3},S_{4}$ arising there one obtains, in our case, that
$$ S_{1}, S_{2}', S_{3} \ll \sqrt{r}x^{1/10}\log^{2}x \frac{x}{dy} \ll x^{9/20}\log^{2}x , $$
$$ \sum_{2^{l}x^{\eta} \leq r \leq 2^{l+1}x^{\eta}} \frac{1}{\phi(r)} \sum_{\substack{\chi^{*} \; (\textrm{mod } r), \\ \chi^{*} \; \textrm{primitive}}} S_{2}'' \ll \left(2^{l}x^{\eta} + \frac{\sqrt{x/d}}{x^{1/20}} + x^{3/20} + \frac{\sqrt{x/d}}{2^{l}x^{\eta}} \right)\sqrt{x/d}\log^{9/2}x , $$
$$ \sum_{2^{l}x^{\eta} \leq r \leq 2^{l+1}x^{\eta}} \frac{1}{\phi(r)} \sum_{\substack{\chi^{*} \; (\textrm{mod } r), \\ \chi^{*} \; \textrm{primitive}}} S_{4} \ll \left(2^{l}x^{\eta} + \frac{\sqrt{x/d}}{x^{1/20}} + \frac{\sqrt{x/d}}{2^{l}x^{\eta}} \right)\sqrt{x/d}\log^{9/2}x . $$
Here we recall that $r \ll \sqrt{x}$ in Theorem 1, and $y > x^{9/10}$ in this subsection. Remembering that we must still multiply by a factor $O(1/\log x)$ that arose from partial summation, these estimates all suffice to give the bound (\ref{largeylargesieve}).
\begin{flushright}
Q.E.D.
\end{flushright}

\appendix
\section{The sums in the exponents}

\subsection{Proof of Lemma 1}
First we make some observations that will make the main part of the proof run more smoothly. We may certainly assume that $\epsilon \geq 1/\log z$, because if it isn't then the bound in Lemma 1 is trivial. We may also assume that $\chi$ is a primitive Dirichlet character, because otherwise we can replace it by the primitive character it is induced from, at the cost of an error term that is
$$ \ll \sum_{p \mid q} \frac{\log p}{p^{\sigma}} \ll \sum_{p \leq 10\log q} \frac{\log p}{p^{\sigma}} \ll \sum_{p \leq 10\log q} \frac{\log p}{p^{0.1}} \ll \log^{0.9}q \;\;\; \textrm{ if } 0.1 \leq \sigma \leq 1, $$
and is
$$ \ll \sum_{p \mid q} \log p \left[\frac{\log z}{\log p}\right] \ll \frac{\log z \log q}{\log\log(q+2)} \;\;\; \textrm{ if } 0 \leq \sigma \leq 0.1, \textrm{ say}. $$
Bearing in mind that we have $z \geq (Hr)^{C} \geq H^{5}$ and $\epsilon \leq 1/2$ in Lemma 1, the second of these terms is
$$ \ll H\log^{2}z + \frac{\log^{2}q}{H} \ll z^{1/4} + \frac{\log^{2}(qzH)}{H} \ll \frac{z^{1-\sigma-0.9\epsilon}}{1-\sigma} + \frac{z^{1-\sigma} \log^{2}(qzH)}{(1-\sigma)H}, $$
so in any case the error term may be absorbed into the right hand side of Lemma 1. Finally, we may assume that $\sigma+it$ is not a zero of $L(s,\chi)$, because if it is then we can replace $\sigma+it$ by an arbitrarily close point that is not a zero, which will have a negligible effect on the left hand side in the statement of the lemma.

\vspace{12pt}
Now a classical explicit formula, reproduced as e.g. Theorem 12.10 of Montgomery and Vaughan~\cite{mv}, implies that if $z,T \geq 2$ and if $\chi$ is a primitive non-principal Dirichlet character then
$$ \sum_{n \leq z} \Lambda(n)\chi(n) = - \sum_{\rho, \atop |\Im(\rho)| \leq T} \frac{z^{\rho}}{\rho} + C(\chi) + O(\log z) + O(\frac{z \log^{2}(rzT)}{T}), $$
where $r$ is the conductor of $\chi$ and
$$ C(\chi) := \frac{L'(1,\overline{\chi})}{L(1,\overline{\chi})} + \log(r/2\pi) - \gamma , $$
with $\gamma$ denoting Euler's constant. The proof of this formula can be modified in a straightforward way to show that, when $0 \leq \sigma < 1$ and $\sigma + it \neq 0$ is not a zero of $L(s,\chi)$,
\begin{eqnarray}
\sum_{n \leq z} \frac{\Lambda(n)\chi(n)}{n^{\sigma + it}} & = & - \sum_{\rho, \atop |\Im(\rho)-t| \leq T} \frac{z^{\rho-\sigma-it}}{\rho - \sigma - it} + (1-a(\chi))\frac{z^{-\sigma-it}}{\sigma + it} -\frac{L'(\sigma+it,\chi)}{L(\sigma+it,\chi)} + O(\frac{\log z}{z^{\sigma}}) \nonumber \\
&& + O(\frac{z^{1-\sigma} \log^{2}(rzT(|t|+1))}{T}), \nonumber
\end{eqnarray}
where $a(\chi)$ is zero or one according as $\chi(-1)$ is 1 or $-1$. (Here the term $(1-a(\chi))z^{-\sigma-it}/(\sigma+it)$ arises because, if $\chi(-1)=1$, the function $L(s,\chi)$ has a zero at $s=0$.)

\vspace{12pt}
At this point we shall divide the proof of Lemma 1 into two cases, according to the relative sizes of $1-\sigma$ and of $\epsilon$ (the width of the hypothesised zero-free region):
\begin{enumerate}
\item if $1-\sigma \leq 0.99\epsilon$;

\item if $1-\sigma > 0.99\epsilon$.
\end{enumerate}

\vspace{12pt}
In the first case, if we choose $T=H/2$ in the preceding discussion, and note that we have $0.505 \leq 1-0.99\epsilon \leq \sigma < 1$ and $|t| \leq H/2$ in Lemma 1, we find that
$$ \left|\sum_{n \leq z} \frac{\Lambda(n)\chi(n)}{n^{\sigma + it}}\right| \ll \sum_{\rho, \atop |\Im(\rho)| \leq H} \frac{z^{\Re(\rho)-\sigma}}{|\rho - \sigma - it|} + \left|\frac{L'(\sigma+it,\chi)}{L(\sigma+it,\chi)}\right| + \frac{\log z}{z^{\sigma}} + \frac{z^{1-\sigma} \log^{2}(rzH)}{H}. $$
Next, a direct modification of the proof of Lemma 3 in the author's paper~\cite{harper3} (replacing $1+1/\log q$ there by $1+\epsilon$, and breaking the sums over zeros according as $|\Im(\rho)| \leq H$, rather than $|\Im(\rho)| \leq q$) shows that
$$ \frac{L'(\sigma+it,\chi)}{L(\sigma+it,\chi)} = O(1/\epsilon + \log(rH)). $$
Keeping in mind that, by assumption, every term $\rho$ in the sum satisfies $\Re(\rho) \leq 1 - \epsilon$, we also have
$$ \sum_{\rho, \atop |\Im(\rho)| \leq H} \frac{z^{\Re(\rho)-\sigma}}{|\rho - \sigma - it|} \ll z^{1/2-\sigma} \sum_{\rho : \Re(\rho) \leq 1/2, \atop |\Im(\rho)| \leq H} \frac{1}{1+|\rho-it|} + \frac{z^{1-\sigma}}{\epsilon} \sum_{k=1}^{[1/2\epsilon]} z^{-k\epsilon} \sum_{\rho : \Re(\rho) > 1-(k+1)\epsilon, \atop |\Im(\rho)| \leq H} 1, $$
since $|\rho-\sigma-it| \geq \max\{|\Re(\rho)-\sigma|,|\Im(\rho)-t|\} \gg \max\{\epsilon,|\Im(\rho)-t|\}$ in this case. Now standard results on the vertical distribution of zeros of $L(s,\chi)$, as in e.g. Theorem 10.17 of Montgomery and Vaughan~\cite{mv}, show that the first sum is $O(\log^{2}(rH))$. Moreover, the log-free zero-density estimate in Zeros Result 1 shows the second sum is
$$ \ll \sum_{k=1}^{[1/2\epsilon]} z^{-k\epsilon} (rH)^{3(k+1)\epsilon} \ll \sum_{k=1}^{[1/2\epsilon]} z^{-0.9k\epsilon} \ll z^{-0.9\epsilon}, $$
provided the value of $C > 0$ in Lemma 1 (for which $z \geq (rH)^{C}$) was chosen large enough. Here we used our assumption that $\epsilon \geq 1/\log z$ to sum the geometric progression. Putting all of this together, and remembering that we have $\epsilon \gg 1-\sigma$ in this first case, we see
$$ \left|\sum_{n \leq z} \frac{\Lambda(n)\chi(n)}{n^{\sigma + it}}\right| \ll \frac{z^{1-\sigma-0.9\epsilon}}{1-\sigma} + \log(rH) + \frac{1}{\epsilon} + \frac{z^{1-\sigma} \log^{2}(rzH)}{H}, $$
which suffices for the bound claimed in Lemma 1.

\vspace{12pt}
In the second case of the proof, where $1-\sigma > 0.99\epsilon$, we shall take a slightly more ``low-tech'' approach. Thus we have
\begin{eqnarray}
\left|\sum_{n \leq z} \frac{\Lambda(n)\chi(n)}{n^{\sigma + it}}\right| & \leq & \sum_{n \leq z^{1/100}} \frac{\Lambda(n)}{n^{\sigma}} + \sum_{j=0}^{[99\log z /(100\log 2)]} \left|\sum_{2^{j}z^{1/100} < n \leq \min\{2^{j+1}z^{1/100},z\}} \frac{\Lambda(n) \chi(n)}{n^{\sigma+it}} \right| \nonumber \\
& \ll & \frac{z^{(1-\sigma)/100}}{1-\sigma} + \sum_{j=0}^{[99\log z /(100\log 2)]} \frac{1}{(2^{j}z^{1/100})^{\sigma}} \max_{m \leq 2^{j+1}z^{1/100}} \left|\sum_{2^{j}z^{1/100} < n \leq m} \frac{\Lambda(n) \chi(n)}{n^{it}} \right|, \nonumber
\end{eqnarray}
where the first line is simply the triangle inequality, and the second line uses Abel's partial summation lemma. Note that
$$ z^{(1-\sigma)/100} = z^{1-\sigma-0.99(1-\sigma)} \leq z^{1-\sigma-0.9\epsilon} $$
in this case, which is acceptable for Lemma 1. We will show that, under the hypotheses of Lemma 1, each subsum in the sum over $j$ is $\ll (2^{j+1}z^{1/100})^{1-0.9\epsilon} + 2^{j+1}z^{1/100}\log^{2}(rzH)/H$, which the reader may check is sufficient to establish the bound claimed in the lemma.

In fact we have already done almost all of the necessary work. The explicit formula that we stated above implies that, for any $X \geq 2$ and any $0 < |t| \leq H/2$,
$$ \left|\sum_{n \leq X} \frac{\Lambda(n)\chi(n)}{n^{it}} \right| \ll \sum_{\rho, \atop |\Im(\rho)| \leq H} \frac{X^{\Re(\rho)}}{|\rho - it|} + \left|(1-a(\chi))\frac{X^{-it}}{it} -\frac{L'(it,\chi)}{L(it,\chi)} \right| + \log X + \frac{X\log^{2}(rXH)}{H}. $$
Moreover, exploiting the functional equation for $L(s,\chi)$, using e.g. formulae (12.9) and (C.17) of Montgomery and Vaughan~\cite{mv}, we find
\begin{eqnarray}
-\frac{L'(it,\chi)}{L(it,\chi)} & = & \frac{L'(1-it,\overline{\chi})}{L(1-it,\overline{\chi})} + \log(r/2\pi) + \frac{\Gamma'(1-it)}{\Gamma(1-it)} - \frac{\pi}{2}\cot((\pi/2)(it+a(\chi))) \nonumber \\
& = & \frac{L'(1-it,\overline{\chi})}{L(1-it,\overline{\chi})} + O(\log(r(|t|+1))) -\frac{1}{it+a(\chi)} + O(1), \nonumber
\end{eqnarray}
and so for $0 < |t| \leq H/2$ we have
\begin{eqnarray}\label{appexplicit}
\left|\sum_{n \leq X} \frac{\Lambda(n)\chi(n)}{n^{it}} \right| \ll \sum_{\rho, \atop |\Im(\rho)| \leq H} \frac{X^{\Re(\rho)}}{|\rho - it|} + \left|\frac{L'(1-it,\overline{\chi})}{L(1-it,\overline{\chi})} \right| + \log(rXH) + \frac{X\log^{2}(rXH)}{H}
\end{eqnarray}
This also holds when $t=0$, that being the standard case that we quoted at the very beginning of this section. In addition, the zero-free region hypothesised in Lemma 1 implies that any ``exceptional'' real zero of $L(s,\overline{\chi})$ is $\leq 1-\epsilon \leq 1-1/\log z$, so standard results (as in Theorem 11.4 of Montgomery and Vaughan~\cite{mv}, for example) imply that $|L'(1-it,\overline{\chi})/L(1-it,\overline{\chi})| \ll \epsilon^{-1} + \log(r(|t|+1)) \ll \log(rzH)$.

Finally, if $z^{1/100} \leq X \leq z$ then, as we did earlier, we can use the bound (\ref{appexplicit}) and the log-free zero-density estimate from Zeros Result 1 to conclude that
$$ \left|\sum_{n \leq X} \frac{\Lambda(n)\chi(n)}{n^{it}} \right| \ll \sqrt{X}\log^{2}(rzH) + X^{1-0.9\epsilon} + \log(rzH) + \frac{X\log^{2}(rzH)}{H}, $$
provided the constant $C > 0$ in Lemma 1 was chosen large enough. The first three terms here are all $\ll X^{1-0.9\epsilon}$, (bearing in mind that $0 < \epsilon \leq 1/2$ and $z \geq (rH)^{C}$), and applying this estimate for $X=2^{j}z^{1/100}$ gives the bound we wanted.
\begin{flushright}
Q.E.D.
\end{flushright}

\subsection{Proof of Lemma 2}
We follow the proof of Lemma 1 closely, with only two changes. Firstly, when we argued that we could replace $\chi$ by the primitive character it is induced from, we required the assumption that $z \geq (Hr)^{C}$ when $0 \leq \sigma \leq 0.1$, say. However, now we can argue that the error term arising there is
$$ \ll \frac{\log z \log q}{\log\log(q+2)} \ll \frac{z^{1-\sigma-0.95\epsilon}\log^{2}(qzH)}{1-\sigma} , $$
which is acceptable for Lemma 2. Secondly, we shall give simpler treatments of some of the sums over zeros in the proof to replace the appeal to a log-free zero-density estimate, which we cannot use successfully having dropped the assumption that $z \geq (Hr)^{C}$.

In the first case of the proof, where $1-\sigma \leq 0.99\epsilon$, we note that
$$ \sum_{\rho, \atop |\Im(\rho)| \leq H} \frac{z^{\Re(\rho)-\sigma}}{|\rho - \sigma - it|} \ll \frac{z^{1-\epsilon-\sigma}}{\epsilon} \sum_{\rho, \atop |\Im(\rho)| \leq H} \frac{1}{1+|\rho - it|} \ll \frac{z^{1-\epsilon-\sigma}}{1-\sigma} \log^{2}rH , $$
in view of standard results on the vertical distribution of zeros. This suffices for the bound claimed in the lemma.

In the second case of the proof, where $1-\sigma > 0.99\epsilon$, it suffices to show that
$$ \sum_{\rho, \atop |\Im(\rho)| \leq H} \frac{X^{\Re(\rho)}}{|\rho - it|} \ll X^{1-\epsilon} \log^{2}(rzH) $$
when $z^{1/100} \leq X \leq z$ (say), and then apply this estimate with $X=2^{j}z^{1/100}$ as in the proof of Lemma 1. However, we immediately see that
$$ \sum_{\rho, \atop |\Im(\rho)| \leq H} \frac{X^{\Re(\rho)}}{|\rho - it|} \ll \sqrt{X}\log^{2}(rzH) + X^{1-\epsilon}\sum_{\rho : \Re(\rho) > 1/2, \atop |\Im(\rho)| \leq H} \frac{1}{1+|\rho - it|} \ll X^{1-\epsilon}\log^{2}(rzH) , $$
as required.
\begin{flushright}
Q.E.D.
\end{flushright}

\vspace{12pt}
\noindent {\em Acknowledgements.} The author would like to thank Andrew Granville for his help with the literature on smooth numbers in arithmetic progressions.

\end{document}